\documentclass[a4paper]{article}

\usepackage[margin=1in]{geometry} 
\usepackage[T1]{fontenc}
\usepackage{amsmath,amssymb,mathtools,graphicx,amsthm}
\usepackage{bm,soul}
\usepackage{algorithm}
\usepackage{algorithmic}
\usepackage{dsfont,bbm}
\usepackage{enumitem}
\usepackage{float}  
\usepackage{caption}
\usepackage{cancel}
\usepackage{siunitx}
\usepackage[scr=boondoxo,scrscaled=1.05]{mathalfa}
\usepackage{physics}
\usepackage{url}
\usepackage{xcolor,soul}
\colorlet{darkred}{red!85!black}
\colorlet{darkgreen}{green!50!black}
\colorlet{darkblue}{blue!60!black}
\usepackage[
colorlinks   = true, 
urlcolor     = darkgreen, 
linkcolor    = darkblue, 
citecolor   = darkred 
]{hyperref}
\usepackage[maxnames=100]{biblatex}
\newtheorem{proposition}{Proposition}[section]
\newtheorem*{remark}{Remark}
\newenvironment{idea}{\noindent\textbf{Idea of the proof:}}{\hfill$\square$ \newline}
\newcommand{\tf}{t_{f}}
\newcommand{\ti}{t_{\iota}}

\title{On the numerical integration of the Fokker-Planck equation driven by a mechanical force and the Bismut-Elworthy-Li formula}
\author{Julia Sanders$^1$ \and Paolo Muratore-Ginanneschi$^1$}

\date{
	$^1$\quad Department of Mathematics and Statistics of the University of Helsinki, 00014 Helsinki, Finland \\ 
    Correspondence:  julia.sanders@helsinki.fi
}

\begin{document}
\maketitle
\begin{abstract}
Optimal control theory aims to find an optimal protocol to steer a system between assigned boundary conditions while minimizing a given cost functional in finite time. Equations arising from these types of problems are often non-linear and difficult to solve numerically. In this note, we describe numerical methods of integration for two partial differential equations that commonly arise in optimal control theory: the Fokker-Planck equation driven by a mechanical potential for which we use Girsanov theorem; and the Hamilton-Jacobi-Bellman, or dynamic programming, equation for which we find the gradient of its solution using the Bismut-Elworthy-Li formula. The computation of the gradient is necessary to specify the optimal protocol. Finally, we give an example application of the numerical techniques to solving an optimal control problem without spacial discretization using machine learning.\end{abstract}

\section{Introduction: optimal control in stochastic thermodynamics and the Fokker-Planck equation.}

Thermodynamic transitions at nanoscale occur in highly fluctuating environments. For instance, nanoscale bio-molecular machines operate within a power output range between $10^{-16}\,\mathrm{W}$ to $10^{-17}\,\mathrm{W}$ per molecule while experiencing random environmental buffeting of approximately $10^{-8}\,\mathrm{W}$ at room temperature~\cite{KaLeZe2007}. Nanomachines experience topological randomness as their motion occurs in inherently non-smooth surroundings due to the fact that machine constituent dimensions are close to those of the atom~\cite{KaLeZe2007}. The dynamics of nanosystems, therefore, need to be described in terms of stochastic \cite{KleF2005} or, more generally, random differential equations \cite{ArnL1998}. Consequently, the laws of macroscopic thermodynamics are replaced by identities involving functions of indicators of the state of the system that are naturally expressed by stochastic processes. Addressing fundamental and technological questions of nanoscale physics has thus propelled interest in the field of stochastic thermodynamics over the last years \cite{SekK2010,SeiU2012,PePi2020,FoJaCa2022}.

A class of important questions in stochastic thermodynamics revolve around finding efficient protocols that natural or artificial nanomachines adopt to perform useful work at nanoscale. Optimal control theory provides a natural mathematical formulation for these type of questions \cite{BecJ2021}.  For instance, conversion of chemical energy into mechanical work typically imply steering the system probability distribution between two assigned values. Schr\"odinger bridge problems \cite{SchE1931} (English translation in \cite{ChMGSc2021}) and their extensions, see e.g. \cite{TodE2009b,LeoC2014,ChGePa2021}, depict this idea mathematically. In these types of problems, protocols optimizing a given functional of the stochastic process describing the state of a nanosystem are determined by solving a pair of coupled Hamilton-Jacobi-Bellman equation \cite{FlSo2006} and a Fokker-Planck equation. 
The solution of the Hamilton-Jacobi-Bellman equation determines the value of the optimal action (force) steering the dynamics at any instant of time in the control horizon. However, the boundary condition at the end of the control horizon is assigned on the system probability distribution. Solving the Fokker-Planck equation thus becomes necessary to fully determine optimal protocols.

Possibly the most prominent physical application of such setup is the derivation of a tight lower bound on the entropy production in classical stochastic thermodynamics \cite{AuGaMeMoMG2012}.  Remarkably, when the system dynamics is modeled by Langevin-Smoluchowski (overdamped) dynamics, the problem maps into the Monge-Amp\`ere-Kantorovich equations and becomes essentially integrable \cite{VilC2009}. This allows one to extricate relevant quantitative information
about molecular engine efficiency \cite{ScSe2008,ScSe2008b,MGSc2015,MaRoDiPePaRi2016} and minimum heat release \cite{PrEhBe2020}. 
The situation is, however, more complicated for more realistic models of nanosystem evolution. Specifically, if we adopt an underdamped (Langevin-Kramers) model of the dynamics \cite{PMG2014,MGSc2014}, then even the equation connecting the optimal mechanical potential to the value function solving the dynamic programming equation is not analytically integrable. 
In the Gaussian case, the solution of a Lyapunov equation \cite{HoJo1991} specifies the optimal mechanical potential \cite{SaBaMG2024}. In general, optimal control duality methods \cite{ChCoGrRe2021} and multiscale perturbative expansions \cite{SaBaMG2024} yield lower bounds on the entropy production and approximate expressions of optimal protocols. More detailed quantitative information calls for exact numerical methods. This is particularly challenging, as integration strategies 
must be adapted to take into account the boundary conditions at the end of the control horizon imposed on the system probability density. Hence, the development of accurate and scale-able methods for numerical integration of the Fokker-Planck equation becomes an essential element of optimization algorithms. 

Traditional numerical methods from hydrodynamics, such as the pseudo-spectral method, see e.g. \cite{CeLaMaVe2001}, are certainly accurate,
but require boundary conditions that are periodic in space, and may not suit problems in stochastic thermodynamics. An even more serious limitation is the exponentially fast increase in computational complexity with the degrees of freedom of the problem.
Monte Carlo averages over Lagrangian particles paths, i.e. realizations of the solution of the stochastic differential equation 
associated to the Fokker-Planck, circumvent the curse of dimensionality, see \cite{MaPMG2001,DoMG2022} for example applications to classical and quantum physics. The drawback is, however, that these methods are best suited for computing expectation values of smooth indicators of the stochastic process. They lack accuracy when computing the probability density itself, as this involves averaging over Dirac distributions. These considerations motivate recent works \cite{MaReOp2020, BoVaEi2022}, which use machine learning methods to construct solutions of the Fokker-Planck equation in the system's state space. These approaches consider the associated probability flow equation \cite{VilC2009,MaReOp2020,BoVaEi2022}, which use the score function (or gradient of the log probability density \cite{HyvA2005}) to turn the Fokker-Planck into a mass conservation equation. The score function can be parametrized by, for example, a neural network \cite{BoVaEi2022}, and the probability density can be recovered through a deterministic transport map. 
 
In this note, we propose a Monte Carlo method adapted to the numerical integration of Fokker-Planck equations of diffusion processes driven by a time-dependent mechanical force. Although mathematically non-generic, they are recurrent in applications of stochastic thermodynamics, as they describe the evolution of a system under a mechanical potential, which may vary in time because of a feedback. We encounter this type of equation in generalized Schr\"odinger bridge problems instantiating refinements of the second law of thermodynamics \cite{AuGaMeMoMG2012,PrEhBe2020,SaBaMG2024}. In this context, the Fokker-Planck equation describes the evolution of the optimal distribution of the state throughout the time interval. Integrating this directly offers a challenging problem, particularly when the driving mechanical potential is non-linear or for systems of high dimension. By using the well-known Girsanov change of measure formula \cite{RiQuRiSc2024}, we couch the solution to the Fokker-Planck in terms of a numerical expectation that can be evaluated from sampled trajectories of the dynamics.

In addition, we also take a look at an application of the Bismut-Elworthy-Li formula \cite{BisJ1984,ElLi1994,ZhaY2010} to compute the gradient of the solution to the Hamilton-Jacobi-Bellman equation. This equation determines the value function, which enforces the system dynamics over the control horizon, and is coupled to the Fokker-Planck through its boundary conditions. Direct access to the gradient of the value function is important, since the stationarity condition in control problems often links the optimal control protocol through the gradient of the value function. The Bismut-Elworthy-Li formula is commonly used in finance, for the calculation of the Greeks derivatives \cite{FoLaLeLiTo1999,FoLaLeLi2001,baños2015}. It has also been used in numerical integration of non-linear parabolic partial differential equations \cite{EHuJeKr2021}. We apply the Bismut-Elworthy-Li formula in the underdamped, or degenerate, dynamics \cite{ZhaY2010}, alongside a numerical example. 

Numerical approaches to the Schrödinger bridge are often iterative. For example, \cite{CaHa2022} turns the problem into a pair of Fokker-Planck equations and iteratively integrates them to recompute the boundary conditions via a proximal operator based numerical integration method \cite{CaHa20}. Machine learning techniques have been used to iteratively solve half-bridge problems \cite{vargas,DeBoThHeAr2021}. We bring together the described Monte Carlo methods for the Fokker-Planck and the Hamilton-Jacobi-Bellman equation in a prototype numerical example to solve a Schrödinger bridge minimizing the Kullback-Leibler divergence from a free diffusion. This is done by an iteration between updating the drift, parametrized by a neural network, with the stationarity condition and the value function with the probability density. Using Monte Carlo integration allows us to compute the update steps quickly and without spacial discretisation.

The next Sections are organized as follows. In Section~\ref{sec:od}, we leverage the Girsanov theorem to express the solution of a Fokker-Planck equation as an expectation which can be evaluated numerically. This is complemented by analytical and numerical examples. In Section~\ref{sec:ud}, we extend this setup into the underdamped dynamics, with an accompanying numerical example from a stochastic optimal control model in the underdamped dynamics. Section ~\ref{sec:od_bel} discusses the application of the Bismut-Elworthy-Li formula to non-degenerate dynamics, with an analytic example application in subsection~\ref{sec:an_hjb_od}, and a  numerical example in subsection~\ref{sec:num_hjb_od} of the Hamilton-Jacobi-Bellman equation coupled to the Fokker-Planck from subsection~\ref{sec:num_fpk_od}. In Section~\ref{sec:ud_bel}, we extend the application of the Bismut-Elworthy-Li formula to the degenerate case. The formula is applied analytically in subsection~\ref{sec:an_hjb_ud}, and numerically in subsection~\ref{sec:num_hjb_ud}. Finally, we give an example use case for the derived formulae, by solving an optimal control problem in the overdamped dynamics by machine learning.

\section{Fokker-Planck for a time-dependent mechanical overdamped diffusion process}
\label{sec:od}

We consider the Langevin-Smoluchowski stochastic differential equation
\begin{align}
	\label{od:LS}
		\mathrm{d}\bm{\mathscr{q}}_{t}=-\mu \,(\bm{\partial}U_{t})(\bm{\mathscr{q}}_{t})\,\mathrm{d}t+\sqrt{\frac{2\,\mu}{\beta}}\,\mathrm{d}\bm{\mathscr{w}}_{t}
	\end{align}
where $ \bm{\mathscr{w}}_{t}$ denotes a standard Wiener process \cite{KleF2005,PavG2014}. The diffusion coefficient $\beta^{-1}$ is proportional to the temperature of the environment surrounding the system. The positive constant $\mu$ is the motility with canonical dimensions of a time over a mass. The drift in (\ref{od:LS}) is the gradient of a time-dependent potential
\begin{align}
	U_{t}(\bm{q})\colon [\ti\,, \,\infty)\,\times\,\mathbb{R}^{d}\mapsto\mathbb{R}_{+}
	\label{od:potential}
\end{align}
that we assume to be sufficiently regular and confining. 

\begin{remark}
    Following a well-established convention in stochastic thermodynamics (see e.g. \cite{ChGa2008}), we denote functional dependence upon time i.e. the dynamical parameter with an underscript. Round brackets express dependence upon state coordinates in configuration or phase space of the diffusion.
\end{remark}
The probability density distribution of the solution of (\ref{od:LS})
at any instant of time $t$ satisfies the Fokker-Planck equation
\begin{align}
	\label{od:FP}
			\partial_{t}\mathtt{p}_{t}(\bm{q})-\mu \,\left \langle\,\bm{\partial}_{\bm{q}}\,,\,(\bm{\partial}U_{t})(\bm{q})\mathtt{p}_{t}(\bm{q})\,\right\rangle-\frac{\mu}{\beta}
		\,\bm{\partial}_{\bm{q}}^{2}\mathtt{p}_{t}(\bm{q})=0
\end{align}
whose solution is fully specified by the assignment of an initial datum at time $t=\ti$.
The assumption of a confining potential (\ref{od:potential}) guarantees that the probability density is integrable in $\mathbb{R}^{d}$.
The connection between (\ref{od:LS}) and (\ref{od:FP}) stems from the representation of the transition probability density as a Monte Carlo average:
\begin{align}
	\mathtt{p}_{t}(\bm{q})=\int_{\mathbb{R}^{d}}\mathrm{d}^{d}\bm{y}\,\operatorname{E}_{\mathcal{P}}\left(\delta^{(d)}(\bm{q}-\bm{\mathscr{q}}_{t})\big{|}\bm{\mathscr{q}}_{\ti}=\bm{y}\right)\,\mathtt{p}_{\ti}(\bm{q})
	\label{od:CK}
\end{align}
The expectation value $\operatorname{E}_{\mathcal{P}}$ is over the probability measure $\mathcal{P}$ weighing the realizations of the solutions of (\ref{od:LS}).
The singular nature of the Dirac delta distribution prevents accurate evaluation of the transition probability density as a Monte Carlo average. For this reason, we look for the solution in the form
\begin{align}
	\mathtt{p}_{t}(\bm{q})=e^{-\beta\,U_{t}(\bm{q})}\,\mathtt{f}_{t}(\bm{q})
	\label{od:A}
\end{align}
Upon inserting into (\ref{od:FP}), we arrive at
\begin{align}
	\partial_{t}\mathtt{f}_{t}(\bm{q})+\mu\,\left \langle\,(\bm{\partial}U_{t})(\bm{q})\,,(\bm{\partial}\mathtt{f}_{t})(\bm{q})\,\right\rangle
	-\frac{\mu}{\beta} 	\,\bm{\partial}_{\bm{q}}^{2}\mathtt{f}_{t}(\bm{q})=\beta\,\mathtt{f}_{t}(\bm{q})\,\partial_{t}U_{t}(\bm{q})
	\label{od:eq}
\end{align}
\begin{proposition}
	The solution of (\ref{od:FP}) admits the representation
	\begin{align}
		\mathtt{p}_{t}(\bm{q})=e^{-\beta\,U_{t}(\bm{q})}\,\operatorname{E}_{\mathcal{P}^{\flat}}\left(	\mathtt{p}_{\ti}(\bm{\mathscr{q}}_{\ti})\,e^{\beta\,U_{\ti}(\bm{\mathscr{q}}_{\ti})+\beta\int_{\ti}^{t}\mathrm{d}s\,\partial_{s}U_{s}(\bm{\mathscr{q}}_{s})}
		\Big{|}\bm{\mathscr{q}}_{t}=\bm{q}\right)
		\label{od:MC}
	\end{align}
	where ${\mathcal{P}^{\flat}}$ is the probability measure over the paths of the backward diffusion process
	\begin{align}
    \begin{split}
		&\mathrm{d}^{\flat}\bm{\mathscr{q}}_{t}:=\bm{\mathscr{q}}_{t}-\bm{\mathscr{q}}_{t-\mathrm{d}t}=\mu\,(\bm{\partial}U_{t})(\bm{\mathscr{q}}_{t})\,\mathrm{d}t+\sqrt{\frac{2\,\mu}{\beta}}\,\mathrm{d}^{\flat}\bm{\mathscr{w}}_{t}
        \\
       & \mathrm{d}^{\flat}\bm{\mathscr{w}}_{t}=\bm{\mathscr{w}}_{t}-\bm{\mathscr{w}}_{t-\mathrm{d}t}
        \end{split}
		\label{od:bsde}
	\end{align}
    naturally complemented by conditions assigned for some $\tf\ge t$.
\end{proposition}
\begin{idea} We start by recalling that, for any test function $\texttt{F}_{t}$, It\^o's lemma for backward differentials yields
	\begin{align}
		\mathrm{d}^{\flat}\mathtt{F}_{t}(\bm{\mathscr{q}}_{t})=\mathrm{d}t\left(\partial_{t}\mathtt{F}_{t}(\bm{\mathscr{q}}_{t})
		+\mu \,\left \langle\,(\bm{\partial}U_{t})(\bm{q})\,,(\bm{\partial}\mathtt{F}_{t})(\bm{q})\,\right\rangle
		-\frac{\mu}{\beta} 	\bm{\partial}_{\bm{q}}^{2}\mathtt{F}_{t}(\bm{q})
		\right)+\sqrt{\frac{2\,\mu}{\beta}}\,\left \langle\,\mathrm{d}^{\flat}\bm{\mathscr{w}}_{t}\,,(\bm{\partial}\mathtt{F}_{t})(\bm{q})\,\right\rangle
		\nonumber
	\end{align}
	We emphasize that $\mathrm{d}^{\flat}\bm{\mathscr{w}}_{t}$ is just a standard Wiener process but evolving backward in time. 
    In the stochastic analysis jargon,
    \eqref{od:bsde} is a diffusion process with respect to a backward filtration as in, e.g., \cite{KunH1982}. As well known, the stochastic integrals and martingale properties are the same as in forward calculus once one exchanges the pre-point rule with the post point, see e.g. \cite{MeyP1982}.
    Let us define the auxiliary function
	\begin{align}
		\mathtt{g}_{t}(\bm{q})=e^{\beta\int_{t}^{\tf}\mathrm{d}s\,\partial_{s}U_{s}(\bm{q})}\mathtt{f}_{t}(\bm{q})
		\nonumber
	\end{align}
    Then, It\^o's lemma and (\ref{od:eq}) immediately imply
	\begin{align}
		\mathrm{d}^{\flat}\mathtt{g}_{t}(\bm{\mathscr{q}}_{t})&=-\,e^{\beta\int_{t}^{\tf}\mathrm{d}s\,\partial_{s}U_{s}(\bm{\mathscr{q}}_{s})}\mathtt{f}_{t}(\bm{\mathscr{q}}_{t})\,\beta\,\partial_{t}U_{t}(\bm{\mathscr{q}}_{s})\,\mathrm{d}t
		+e^{\beta\int_{t}^{\tf}\mathrm{d}s\,\partial_{s}U_{s}(\bm{\mathscr{q}}_{s})}\,\mathrm{d}\mathtt{f}_{t}(\bm{\mathscr{q}}_{t})
		\nonumber\\
		&=e^{\beta\int_{t}^{\tf}\mathrm{d}s\,\partial_{s}U_{s}(\bm{\mathscr{q}}_{s})} \sqrt{\frac{2\,\mu}{\beta}}\,\left \langle\,\mathrm{d}^{\flat}\bm{\mathscr{w}}_{t}\,,(\bm{\partial}\mathtt{f}_{t})(\bm{q})\,\right\rangle
		\nonumber
	\end{align}
	The equation tells us that the auxiliary function is a local martingale of the backward diffusion (see e.g. Chapter~7 of \cite{KleF2005}). Since we assume that the confining potential also guarantees integrability, we infer that the martingale property
	\begin{align}
		\operatorname{E}_{\mathcal{P}^{\flat}}\left(\mathtt{g}_{\tf}(\bm{\mathscr{q}}_{\tf})\,\big{|}\,\bm{\mathscr{q}}_{\tf}=\bm{q}\right)=
		\operatorname{E}_{\mathcal{P}^{\flat}}\left(\mathtt{g}_{\ti}(\bm{\mathscr{q}}_{\ti})\,\big{|}\,\bm{\mathscr{q}}_{\tf}=\bm{q}\right)
		\nonumber
	\end{align}
	must also hold true. By construction, we know that
	\begin{align}
		\operatorname{E}_{\mathcal{P}^{\flat}}\left(\mathtt{g}_{\tf}(\bm{\mathscr{q}}_{\tf})\big{|}\bm{\mathscr{q}}_{\tf}=\bm{q}\right)=\mathtt{f}_{\tf}(\bm{q})
		\nonumber
	\end{align}
	and we conclude
	\begin{align}
		\mathtt{f}_{\tf}(\bm{q})=	\operatorname{E}_{\mathcal{P}^{\flat}}\left(\mathtt{g}_{\ti}(\bm{\mathscr{q}}_{\ti})\,\big{|}\,\bm{\mathscr{q}}_{\tf}=\bm{q}\right)=\operatorname{E}_{\mathcal{P}^{\flat}}\left(e^{\beta\int_{\ti}^{\tf}\mathrm{d}s\,\partial_{s}U_{s}(\bm{q})}\mathtt{f}_{\ti}(\bm{\mathscr{q}}_{\ti})\,\big{|}\,\bm{\mathscr{q}}_{\tf}=\bm{q}\right )
		\nonumber
	\end{align}
	Replacing $\tf$ with $t$ in the above chain of identities completes the proof. 
\end{idea}
The upshot is that we can use Feynman-Kac formula over a backward diffusion to compute the solution of a forward Fokker-Planck equation. Next, we take advantage of Girsanov's change of measure formula (see e.g. Chapter~10 of \cite{KleF2005} or 3.5 of \cite{PavG2014}) to evaluate the conditional expectation in (\ref{od:MC}) directly over the paths of the Wiener process, or, more generally, over the paths of any diffusion that generates a measure with respect to which $\mathcal{P}^{\flat}$ is absolutely continuous. Girsanov’s change of measure formula is thus the basis of statistical inference for diffusion processes, see e.g. \cite{SoeM2012}. We emphasize that we make use of Girsanov formula while dealing with backward diffusions as e.g. in \cite{MeyP1982}. As time is evolving from a larger to smaller value, correspondingly, the role of ``past'' and ``future'' events must be exchanged.
\begin{remark}
    As increments of any Wiener process are independent, from now on we write
    \begin{align}
    \mathrm{d}^{\flat}\bm{\mathscr{w}}_{t}=\mathrm{d}\bm{\mathscr{w}}_{t}
    \nonumber
    \end{align}
    to alleviate the notation.
\end{remark}

\subsection{Use of Girsanov's formula}

We denote by $\mathcal{Q}$ the probability measure over the path of
\begin{align}
	\mathrm{d}^{\flat}\bm{\mathscr{q}}_{t}=\sqrt{\frac{2\,\mu}{\beta}}\,\mathrm{d}\bm{\mathscr{w}}_{t}
	\label{od:free}
\end{align}
Our aim is to use Girsanov's formula to express expectations with respect to the path measure $\mathcal{P}^{\flat}$ of (\ref{od:bsde}) in terms of expectations with respect to $\mathcal{Q}$:
\begin{align}
	\mathtt{p}_{t}(\bm{q})=e^{-\beta\,U_{t}(\bm{q})}\,\operatorname{E}_{\mathcal{Q}}\left(\mathtt{p}_{\ti}(\bm{\mathscr{q}}_{\ti})\,e^{\beta U_{\ti}(\bm{\mathscr{q}}_{\ti})+\beta\int_{\ti}^{t}\mathrm{d}s\,\partial_{s}U_{s}(\bm{\mathscr{q}}_{s})} \frac{\mathrm{d}\mathcal{P}^{\flat}}{\mathrm{d} \mathcal{Q}}\,\Big{|}\,\bm{\mathscr{q}}_{t}=\bm{q}\right)
	\nonumber
\end{align}
where 
\begin{align}
	\frac{\mathrm{d}\mathcal{P}^{\flat}}{\mathrm{d} \mathcal{Q}}=\exp\int_{\ti}^{t}\left(\sqrt{\frac{\beta\,\mu}{2}}\left \langle\,\mathrm{d}\bm{\mathscr{w}}_{s}\overset{\circ}{\,,}\,(\bm{\partial}U_{s})(\bm{\mathscr{q}}_{s})\,\right\rangle
	-\frac{\beta\,\mu}{4}\mathrm{d}s \left\|(\bm{\partial}U_{s})(\bm{\mathscr{q}}_{s})\right\|^{2}
	\right)
	\label{od:Girsanov}
\end{align}
is the Radon-Nikodym derivative.
The symbol $\circ$ emphasizes that we define the stochastic integral in (\ref{od:Girsanov}) using the \emph{post-point} prescription:
\begin{align}
	\int_{\ti}^{t}\left \langle\,\mathrm{d}\bm{\mathscr{w}}_{s}\overset{\circ}{\,,}\,(\bm{\partial}U_{s})(\bm{\mathscr{q}}_{s})\,\right\rangle
	:=\lim_{\substack{\mathrm{d}t \downarrow 0 \\ N\mathrm{d}t=t-\ti } }\sum_{i=1}^{N} \left \langle\,\bm{\mathscr{w}}_{t_{i+1}}-\bm{\mathscr{w}}_{t_{i}}\,,(\bm{\partial}U_{t_{i+1}})(\bm{\mathscr{q}}_{t_{}})\,
		\right\rangle
	\nonumber
\end{align}
i.e., the function $(\bm{\partial}U_{t})(\bm{\mathscr{q}}_{t})$ is evaluated at end of each time interval. Accordingly, (\ref{od:Girsanov}) is a martingale with respect to the backward filtration (i.e. the family of $\sigma$-algebras increasing as $\ti$ decreases) to which we associate the probability measure $\mathcal{Q}$. In other, rougher, words, (\ref{od:Girsanov}) is a martingale conditional to events occurring at times larger or equal to the upper bound of integration $t$.

Writing the stochastic integral in the standard pre-point form allows us to simplify the expression of the probability density.
We notice that (\ref{od:free}) trivially implies 
\begin{align}
	\int_{\ti}^{t}\sqrt{\frac{\beta\,\mu}{2}}\left \langle\,\mathrm{d}\bm{\mathscr{w}}_{s}\overset{\circ}{\,,}\,(\bm{\partial}U_{s})(\bm{\mathscr{q}}_{s})\,\right\rangle
		=
	\int_{\ti}^{t}\frac{\beta\,\mu}{2}\left \langle\,\mathrm{d}^{\flat}\bm{\mathscr{q}}_{s}\overset{\circ}{\,,}\,(\bm{\partial}U_{s})(\bm{\mathscr{q}}_{s})\,\right\rangle
	\nonumber
\end{align}
Next, we use the relation between stochastic integrals in the post-point, mid-point (or Stratonovich, denoted by the $\diamond$-symbol), and pre-point
prescriptions:
 \begin{align}
 	\int_{\ti}^{t}\left \langle\,\mathrm{d}\bm{\mathscr{q}}_{s}\overset{\diamond}{\,,}\,(\bm{\partial}U_{s})(\bm{\mathscr{q}}_{s})\,\right\rangle
 	=
 	\int_{\ti}^{t}\frac{
 		\left \langle\,\mathrm{d}\bm{\mathscr{q}}_{s}\overset{\circ}{\,,}\,(\bm{\partial}U_{s})(\bm{\mathscr{q}}_{s})\,\right\rangle
 		+
 		\left \langle\,\mathrm{d}\bm{\mathscr{q}}_{s}\,,(\bm{\partial}U_{s})(\bm{\mathscr{q}}_{s})\,\right\rangle
 	}{2}
 	\nonumber
 \end{align}
Finally, we recall that ordinary differential calculus holds for the stochastic differentials  in Stratonovich form
\begin{align}
		U_{t}(\bm{\mathscr{q}}_{t})-U_{\ti}(\bm{\mathscr{q}}_{\ti})
	=
	\int_{\ti}^{t}\Bigl(\mathrm{d}s \,\partial_{s}U_{s}(\bm{\mathscr{q}}_{s}) +\left \langle\,\mathrm{d}\bm{\mathscr{q}}_{s}\overset{\diamond}{\,,}\,(\bm{\partial}U_{s})(\bm{\mathscr{q}}_{s})\,\right\rangle\Bigr )
	\nonumber
\end{align}
Putting these observations together, we obtain the following representation of the solution for the Fokker-Planck equation (\ref{od:FP})
\begin{align}
	\mathtt{p}_{t}(\bm{q})=\,\operatorname{E}\left(\mathtt{p}_{\ti}(\bm{\mathscr{q}}_{\ti})\,
	e^{-\frac{\beta}{2}\int_{\ti}^{t}\left(\left \langle\,\mathrm{d}\bm{\mathscr{q}}_{s}\,,\,(\bm{\partial}U_{s})(\bm{\mathscr{q}}_{s})\,\right\rangle
		+\frac{\mu}{2}\,\mathrm{d}s\,\left\|(\bm{\partial}U_{s})(\bm{\mathscr{q}}_{s})\right\|^{2}
		\right)}
	 \,\Big{|}\,\bm{\mathscr{q}}_{t}=\bm{q}\right)
	\label{od:main}
\end{align} 
In practice, this means that, to compute the probability density at the configuration space point $\bm{q}$ at time $t\,\geq\,\ti$, we need to average the initial density over solutions of (\ref{od:free}) evolved backward in time to $\ti$ and weighted by a path-dependent change-of-measure factor.

\subsection{Examples and path integral representation}

Let's summarize the meaning of \eqref{od:main} in words. Formula \eqref{od:MC} tells us that a Fokker-Planck equation of a forward diffusion process with gradient drift, i.e. of the form \eqref{od:FP}, admits a Feynman-Kac representation in terms of a backward diffusion process. This is because we can use the potential
specifying the drift to turn the forward Fokker-Planck into a non-homogeneous backward Kolmogorov equation with respect to the backward diffusion process. This latter equation, as well known, generically specifies a problem with initial data. We now turn to illustrate this fact with two examples.

\subsubsection{Analytical example}

Consider a quadratic potential
\begin{align}
	U_{t}(\bm{q})=\frac{1}{2}\left \langle\,\bm{q}\,,\mathsf{U}_{t}\bm{q}\,\right\rangle
	\nonumber
\end{align}
with $\mathsf{U}_{s} $ a $d\,\times\,d$ real symmetric time dependent matrix. The backward stochastic differential equation \eqref{od:bsde} reduces to
\begin{align}
	\mathrm{d}^{\flat}\bm{\mathscr{q}}_{t}=\mu\,\mathsf{U}_{t}\bm{\mathscr{q}}_{t}\mathrm{d}t+\sqrt{\frac{2\,\mu}{\beta}}\mathrm{d}\bm{\mathscr{w}}_{t}
	\nonumber
\end{align} 
Letting $\mathsf{F}$ denote  the flow solution of the deterministic ordinary differential equation
\begin{align}
\begin{split}
    \frac{\mathrm{d}}{\mathrm{d}t }\mathsf{F}_{t,s}=-\mu\,\mathsf{U}_{t}\,\mathsf{F}_{t,s}
\end{split}    
	\nonumber
\end{align}
then the solution of the backward stochastic differential equation is
\begin{align}
\nonumber
	&\bm{\mathscr{q}}_{t}=\mathsf{F}_{\tf,t}^{\top}\,\bm{q}-\sqrt{\frac{2\,\mu}{\beta}}\int_{t}^{\tf}\mathsf{F}_{s,t}^{\top}\,\mathrm{d}\bm{w}_{s}
	\\
    &\bm{\mathscr{q}}_{\tf}=\bm{q}
 \nonumber
\end{align}
The symbol $\top$ as usual denotes matrix transposition. The corresponding transition probability density is Gaussian with mean
\begin{align}
&\operatorname{E}	\left(\bm{\mathscr{q}}_{t}\big{|}\bm{\mathscr{q}}_{u}=\bm{q}\right )=\mathsf{F}_{u,t}^{\top}\bm{q}&& u\,\geq\,t
	\nonumber
\end{align}
(recalling that, for standard backward differential equations, the martingale property arises upon conditioning on future events \cite{KunH1982}), and variance matrix
\begin{align}
&\operatorname{E}	\left((\bm{\mathscr{q}}_{t}-\operatorname{E}\bm{\mathscr{q}}_{t})\otimes(\bm{\mathscr{q}}_{t}-\operatorname{E}\bm{\mathscr{q}}_{t})\big{|}\bm{\mathscr{q}}_{u}=\bm{q}\right )=
	\frac{2\,\mu}{\beta}\int_{t}^{u}\mathrm{d}s\,\mathsf{F}_{s,t}^{\top}\mathsf{F}_{s,t}&& u\,\geq\,t
	\nonumber
\end{align}
We need to compute
\begin{align}
	f_{t}(\bm{q})=\operatorname{E}\left(f_{\ti}(\bm{q})e^{\beta\int_{\ti}^{t}\mathrm{d}s \frac{\left \langle\,\bm{\mathscr{q}}_{s}\,,\dot{\mathsf{U}}_{s}\bm{\mathscr{q}}_{s}\,\right\rangle}{2}}\big{|}\bm{q}_{t}=\bm{q}\right)
	\nonumber
\end{align}
If we couch this expression into the form of a path integral \cite{LaRoTi1982}, we get
\begin{align}
	f_{t}(\bm{q})=\int_{\bm{\mathscr{q}}_{t}=\bm{q}}\mathcal{D}[\mathsf{q}_{t:\ti}] e^{-\int_{\ti}^{t}\mathrm{d}s\left(\beta \frac{\|\dot{\mathsf{q}}_{s}-\mu\,\mathsf{U}_{s}\bm{\mathsf{q}}_{s}\|^{2}}{4\,\mu}-\beta \frac{\left \langle\bm{\mathsf{q}}_{s}\,,\dot{\mathsf{U}}_{s}\bm{\mathsf{q}}_{s}\,\right\rangle\,}{2}\right)}f(\bm{\mathsf{q}}_{\ti})
	\nonumber
\end{align}
Here $\mathcal{D}[\mathsf{q}_{t:\ti}]$ denotes the limit over finite dimensional approximations over time lattices in $[\ti,t]$ of paths satisfying the terminal condition $\bm{\mathscr{q}}_{t}=\bm{q}$.
We are free to interpret the path integral in the mid-point sense, because any change of discretization generates a path-independent Jacobian 
that can be reabsorbed into the normalization constant. As the integral is Gaussian, we can compute it  by infinite dimensional stationary phase using ordinary differential calculus. We are left with
\begin{align}
	f_{t}(\bm{q})=e^{U_{t}(\bm{q})}\int_{\bm{\mathscr{q}}_{t}=\bm{q}}\mathcal{D}[\mathsf{q}_{t:\ti}] e^{-\int_{\ti}^{t}\mathrm{d}s\left(\beta \frac{\|\dot{\mathsf{q}}_{s}\|^{2}}{4\,\mu}+\beta \frac{2\left \langle\bm{\dot{\mathsf{q}}}_{s}\,,\mathsf{U}_{s}\bm{\mathsf{q}}_{s}\,\right\rangle+\mu\left \langle\mathsf{U}_{s}\bm{\mathsf{q}}_{s}\,,\mathsf{U}_{s}\bm{\mathsf{q}}_{s}\,\right\rangle}{4}\right)}p_{\ti}(\bm{\mathsf{q}}_{\ti})
	\nonumber
\end{align}
We now readily recognize that 
\begin{align}
	\mathtt{p}_{t,s}(\bm{q}|\bm{y})=\int_{\bm{\mathscr{q}}_{s}=\bm{y}}^{\bm{\mathsf{q}}_{t}=\bm{q}}\mathcal{D}[\mathsf{q}_{t:\ti}] e^{-\int_{\ti}^{t}\mathrm{d}s\left(\beta \frac{\|\dot{\mathsf{q}}_{s}\|^{2}}{4\,\mu}+\beta \frac{2\left \langle\bm{\dot{\mathsf{q}}}_{s}\,,\mathsf{U}_{s}\bm{\mathsf{q}}_{s}\,\right\rangle+\mu\left \langle\mathsf{U}_{s}\bm{\mathsf{q}}_{s}\,,\mathsf{U}_{s}\bm{\mathsf{q}}_{s}\,\right\rangle}{4}\right)}
	\nonumber
\end{align}
is the path integral expression of the transition probability density of the forward stochastic differential equation
\begin{align}
	\mathrm{d}\bm{\mathscr{q}}_{t}=-\mu \,\mathsf{U}_{t}\bm{\mathscr{q}}_{t}\mathrm{d}t+\sqrt{\frac{2\,\mu}{\beta}}\mathrm{d}\bm{\mathscr{w}}_{t}
	\nonumber
\end{align}
We therefore recover the Chapman-Kolmogorov representation of the solution of \eqref{od:FP}
\begin{align}
	\mathtt{p}_{t}(\bm{q})=\int_{\mathbb{R}}\mathrm{d}\bm{y}\,\mathtt{p}_{t,s}(\bm{q}|\bm{y})\,	\mathtt{p}_{\ti}(\bm{q})
	\label{od:CKex}
\end{align}
as expected.

\subsubsection{Path-integral representation in general}

The path integral representation of \eqref{od:MC} is
\begin{align}
	\mathtt{p}_{t}(\bm{q})=\int_{\mathbb{R}^{d}}\mathrm{d}^{d}\bm{y}\int_{\bm{\mathscr{q}}_{\ti}=\bm{y}}^{\bm{\mathscr{q}}_{t}=\bm{q}}\mathcal{D}[\mathsf{q}_{t:\ti}] e^{-\int_{\ti}^{t}\mathrm{d}s\left(\beta \frac{\|\dot{\mathsf{q}}_{s}\|^{2}}{4\,\mu}+\beta \frac{\left \langle\bm{\dot{\mathsf{q}}}_{s}\,,(\bm{\partial}U_{s})(\bm{\mathsf{q}}_{s})\,\right\rangle+\mu\|(\bm{\partial}U_{s})(\bm{\mathsf{q}}_{s})\|^{2}}{4}\right)}p_{\ti}(\bm{y})
	\nonumber
\end{align} 
As the stochastic integral term is evaluated in the pre-point representation, the path integral exactly recovers the path integral expression of the transition
probability density. We have thus verified that \eqref{od:CKex} holds in general. We refer the reader unfamiliar with path integral calculus to e.g. \cite{LaRoTi1982}. 

\subsubsection{Numerical Example: Time-Independent Drift}
\begin{figure}
     \centering
    \includegraphics[width=\linewidth]{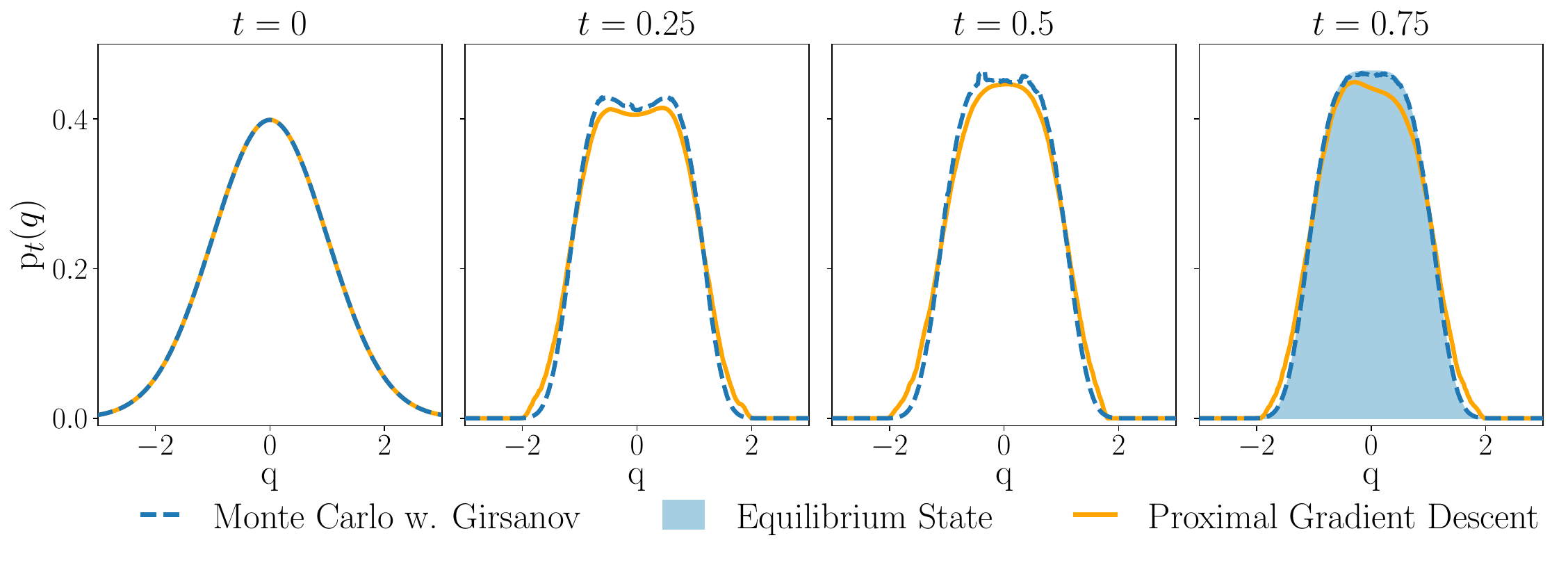}
     \caption{ \label{fig:timeindep_fpk}Solution of a Fokker-Planck equation driven by a mechanical potential~\eqref{ud:FP} computed using Monte Carlo integration via Girsanov formula (dashed blue line). We use $\partial_q U(q) = 2\,q^3$.
     The initial condition is $\mathrm{p}_{\ti}(q)=\frac{1}{\sqrt{2\pi}}\exp(-q^2/2)$ at $\ti=0$. The Girsanov method is compared with an implementation of the "proximal gradient descent" method described in \cite{CaHa20}, shown in orange. For the proximal gradient descent, we use $10^{4}$ samples from the initial distribution and $\gamma=0.05$ as the regularisation parameter, see \cite{CaHa20}. Both methods simulate trajectories of the auxilliary stochastic process~\eqref{od:main} by the Euler-Maruyama scheme with step size $h=10^{-3}$. For the Girsanov theorem approach, we evolve $10^{3}$ trajectories from $10^{4}$ initial points in the interval $[-6,6]$. Resulting distributions are smoothed by convolution with a box filter. We use $\ti=0$ and $\mu=\beta=1$. The expected equilibrium state of the distribution is shown by the shaded area in the final panel at $t=0.75$. In our implementation, the Monte Carlo method of integration is roughly three orders of magnitude faster than the proximal gradient descent. Accompanying code for all figures can be found in the link in the Data Availability statement.}
 \end{figure}

In this section, we demonstrate how the method described can be applied to find a numerical solution of a Fokker-Planck equation driven by a mechanical potential. We consider the Fokker-Planck of the form \eqref{od:FP} with a time-independent drift. By applying the Girsanov theorem as above, we couch the solution of the Fokker-Planck at time $t\in [\ti,\tf]$ into a numerical average of simulated trajectories of the auxiliary dynamics, given by
\begin{equation}\label{num:fpk_expectation}
    \begin{split}
    \mathtt{p}_{t}(\bm{q}) &= \operatorname{E}_{\mathcal{Q}}\left(P_{\iota}(\bm{\mathscr{q}}_{\ti})\exp\biggl(- \int_{\ti}^t\sqrt{\dfrac{\beta\,\mu}{2}}\Bigl\langle\dd \bm{\mathscr{w}}_r\,, \, \bm{\partial}_{\bm{q}}U_{r}(\bm{\mathscr{q}}_{r}) \Bigr\rangle + \dd r \,\dfrac{\beta\,\mu}{4}\Bigl\| \bm{\partial}_{\bm{q}} U_{r}(\bm{\mathscr{q}}_{r})  \Bigr\|^2 \biggr)\,\Big|\, \bm{\mathscr{q}}_{t} = \bm{q}\right) \\
    \end{split}
\end{equation}
where $\mathcal{Q}$ is the measure generated by the backwards diffusion~\eqref{od:free}. We approximate the expectation value numerically by repeated sampling of trajectories of the process~\eqref{od:free}. The trajectories are approximated on a discretization of the time interval $[\ti,\tf]$ given by  
\begin{equation}
    \ti=t_0<t_1<\ldots<t_N=\tf
\end{equation}
Trajectories of \eqref{od:free} are sampled using the Euler-Maruyama scheme
\begin{equation}
    \bm{\mathscr{q}}_{t_{n-1}} = \bm{\mathscr{q}}_{t_n} - \sqrt{\dfrac{2\,\mu\,}{\beta}|t_{n-1}-t_n|}\,\epsilon
\end{equation}
where $\epsilon$ is an increment of Brownian noise, sampled independently from a standard normal distribution. The Girsanov factor $g$ is computed as a running cost
\[g = \sum_{n=0}^{N-1}\dfrac{\beta\,\mu}{4}  |t_{n-1}-t_n| \|(\partial \bm{U})_{t_{n-1-i}} (\bm{\mathscr{q}}_{t_{n-1-i}})|^2 +\sqrt{\dfrac{\mu\, \beta}{2}|t_{n-1}-t_n|} \,(\partial \bm{U}_{t_{n-1-i}}) (\bm{\mathscr{q}}_{t_{n-1-i}})\ \epsilon\]
This computation is summarized in Algorithm~\ref{alg:fpkgirsanov}.  In Fig.~\ref{fig:timeindep_fpk}, we integrate an example Fokker-Planck equation driven by a time-independent mechanical potential in two ways. The results of Algorithm~\ref{alg:fpkgirsanov} are compared to the proximal gradient descent method of \cite{CaHa20}. In this method, the solution is found via gradient descent on the space of probability distributions by solving a proximal fixed point recursion at each time step. Both methods discretize the time interval, but do not require spacial discretization. In our implementation, the Monte Carlo method performs significantly faster. 

\begin{algorithm}
\caption{Integrating Fokker-Planck equation using Girsanov theorem}
\label{alg:fpkgirsanov}
\begin{algorithmic} 
\STATE Initialize $\bm{\mathscr{q}}_{t_n} = \bm{q}\in \mathbb{R}^d$
\STATE Initialize $g = 0$
\STATE Initialize $\partial U_{t_n}$ for $n\in\{0,1,\ldots,N\}$
\STATE Initialize $\delta_{t} = |t_{n-1}-t_n|$ for $n\in\{0,1,\ldots,N\}$
\FOR{ $i$ in $0,\ldots, n-1$}
\STATE Sample Brownian noise: $\epsilon \sim \mathcal{N}(0,1)$
\STATE Evolve one step of \eqref{od:bsde}: $\bm{\mathscr{q}}_{t_{n-i-1}} = \bm{\mathscr{q}}_{t_{n-i}} - \sqrt{2\,\delta_t\mu/\beta }\ \epsilon$ 
\STATE Add to running total: \[g = g +\delta_t \frac{\mu\beta}{4}\,\|(\partial \bm{U}_{t_{n-1-i}})(\bm{\mathscr{q}}_{t_{n-1-i}})\|^2 +\sqrt{\delta_t\frac{\mu\beta}{2}}\,(\partial \bm{U}_{t_{n-1-i}})(\bm{\mathscr{q}}_{t_{n-1-i}})\ \epsilon\]
\ENDFOR
\STATE Return $\mathtt{p}_{t_n}(\bm{q}) = P_{\iota}(\bm{\mathscr{q}}_{\ti})\,e^{-g}$
\end{algorithmic}
\end{algorithm}
 
\subsubsection{Numerical example: "F\"ollmer's drift"}\label{sec:num_fpk_od}
In this section, we apply \eqref{od:MC} to a non-trivial example of gradient drift. Specifically, we consider the F\"ollmer-drift solution of the dynamic Schr\"odinger bridge that steers the system between assigned boundary conditions while minimizing the Kullback-Leibler divergence from a free diffusion~\cite{ChMGSc2021}, given by
\begin{equation}
    \label{eq:kl_divergence}
D_{\mathcal{KL}} = \dfrac{\beta\,\mu}{4}\,\operatorname{E}\left(\int_{\ti}^{\tf} \dd t \,\|(\partial U_t)(\bm{\mathscr{q}}_{t})\|^2 \right)\end{equation}
in a finite time interval $t\in [\ti,\tf]$. The boundary conditions are assigned on the initial and final distributions of the position, denoted $P_{\iota}$ for the initial at time $t=\ti$ and $P_f$ for the final at time $t=\tf$. We consider boundary conditions of the form
\begin{subequations}
\label{eq:boundary_conditions}
\begin{align}\label{eq:boundary_conditions_initial}
P_{\iota}(\bm{q})&=\dfrac{e^{-\beta\, U_{\iota}(\bm{q}))}}{\int_{\mathbb{R}^d}\dd^d \bm{y}\, e^{-\beta\, U_{\iota}(\bm{y}))}}\\
\label{eq:boundary_conditions_final}P_f(\bm{q})&=\dfrac{e^{-\beta\, U_{f}(\bm{q})}}{\int_{\mathbb{R}^d}\dd^d \bm{y}\, e^{-\beta \,U_{f}(\bm{y}))}} \,.  
\end{align}    
\end{subequations}
F\"ollmer-drifts are relevant to machine learning applications see e.g. \cite{TzRa2019,DeBoThHeAr2021}. We refer to \cite{LeoC2014} or \cite{SaBaMG2024} and references therein for further details on mathematics and physics background, respectively.

We summarize how to construct the F\"ollmer drift by solving a Schr\"odinger bridge problem using an iterative method of~\cite{CaHa2022}. In doing so, we also obtain the solution of the Fokker-Planck equation \eqref{od:FP} that we use for comparison with the numerical expression provided by \eqref{od:MC}. The Schr\"odinger bridge problem is formulated as the minimization of a Bismut-Pontryagin functional \cite{SaBaMG2024}. In this framework, we find that the intermediate density $\mathrm{p}$ and a value function $V$ imposing the boundary conditions satisfy the coupled partial differential equations 
\begin{subequations}
\label{eq:coupled_pde_sys}
    \begin{align}      \label{eq:fpk_coupled}
    \partial_{t}\mathtt{p}_{t}(\bm{q})-\mu\left \langle\,\bm{\partial}_{\bm{q}}\,,(\bm{\partial}U_{t})(\bm{q})\mathtt{p}_{t}(\bm{q})\,\right\rangle-\frac{\mu}{\beta}
		\bm{\partial}_{\bm{q}}^{2}\mathtt{p}_{t}(\bm{q})&=0; \\
        \label{eq:hjb_coupled}
    \partial_t V_t(\bm{q}) -\mu\,\left \langle\,(\bm{\partial}U_{t})(\bm{q})\,\,,\bm{\partial}_{\bm{q}} V_t(\bm{q})\right\rangle + \frac{\mu}{\beta} \,\bm{\partial}_{\bm{q}}^{2}  V_t(\bm{q}) + \frac{\beta\,\mu}{4}\Bigl((\bm{\partial}U_{t})(\bm{q})\Bigr)^2 &= 0
    \end{align}
\end{subequations}
along with the stationarity condition 
\begin{equation}
    \label{stationarity_condition}\bm{\partial}_{\bm{q}} V_t(\bm{q}) = \dfrac{2}{\beta}\bm{\partial}_{\bm{q}} U_t(\bm{q})
\end{equation}

We identify \eqref{eq:fpk_coupled} as the Fokker-Planck equation, and \eqref{eq:hjb_coupled} as the Hamilton-Jacobi-Bellman equation which is discussed in later sections. For a known $U$, we can apply Girsanov theorem to integrate  \eqref{eq:fpk_coupled}. We find a reference solution to the system  \eqref{eq:coupled_pde_sys} using an adaptation of the method of~\cite{CaHa2022}, which is briefly described below. The transformation   

\begin{equation}\label{num1:transform}\begin{split}
        \mathtt{p}_{t}(\bm{q})&= \phi_t(\bm{q})\hat\phi_t(\bm{q}) \\
V_t(\bm{q})&=-\log(\phi_t(\bm{q}))\end{split}
\end{equation}
applied to \eqref{eq:coupled_pde_sys} yields the linear coupled equations 
\begin{subequations}
        \begin{align}\label{num:ivpA}
     \partial_t \phi_t(\bm{q}) + \dfrac{\mu}{\beta}\,\bm{\partial}^2_{\bm{q}} \phi_t(\bm{q})&=0 \\ 
     \label{num:ivpB}\partial_t \hat\phi_t(\bm{q}) - \dfrac{\mu}{\beta}\,\bm{\partial}^2_{\bm{q}} \hat\phi_t(\bm{q})&=0 
    \end{align}
\end{subequations}
with boundary conditions 
        \begin{align*}
    \phi_{\tf}(\bm{q}) &= P_f(\bm{q})\,/\, \hat\phi_{\tf}(\bm{q}) \\
    \hat\phi_{\ti}(\bm{q}) &= P_\iota(\bm{q})\,/ \,\phi_{\ti}(\bm{q})
    \end{align*}
We make an initial guess for $\phi_{\ti}$ which we use to integrate \eqref{num:ivpA}, recompute the boundary conditions and then integrate \eqref{num:ivpB}, recomputing $\phi_{\ti}$, and repeating this process until convergence: see~\cite{CaHa2022} or Section 8.2 of~\cite{SaBaMG2024} for a more detailed treatment.   
We reconstruct the value function and intermediate densities using~\eqref{num1:transform}. 

\begin{figure}
     \centering
    \includegraphics[width=\linewidth]{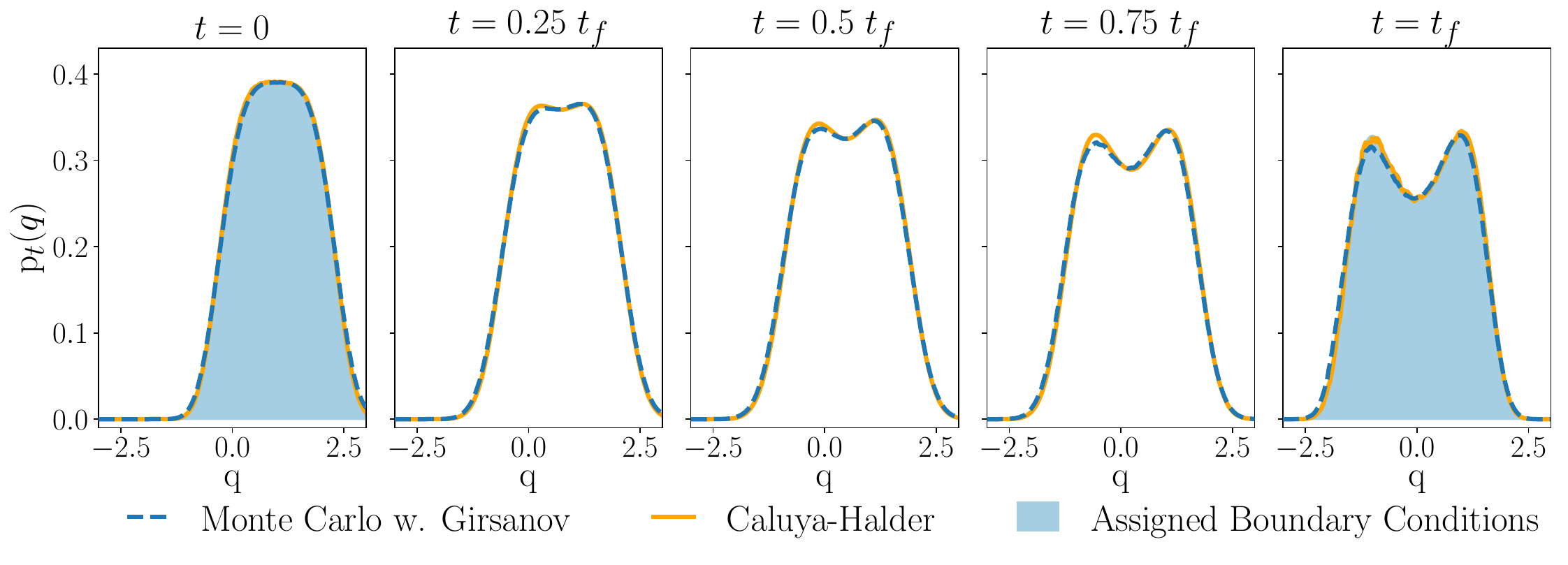}
     \caption{ \label{fig:overdamped_fpk}Solution of a Fokker-Planck driven by a time-dependent mechanical potential computed using Monte Carlo integration via Girsanov formula (dashed blue line). The optimal protocol $U$ and reference solution (orange line) is computed using an iterative method~\cite{CaHa2022}. For the Girsanov theorem approach, we evolve $M = 10\,000$ trajectories from $500$ initial points in the interval $[-3,3]$ with a time step of $h=0.005$ by the Euler-Maruyama scheme. The reference solution uses the iteration method of~\cite{CaHa2022}, where integration of the equations \eqref{num:ivpA} and \eqref{num:ivpB} is also computed as a numerical average of Monte Carlo sampled trajectories, using $5\,000$ initial points from the interval $[-6,6]$. Ten total iterations are performed. Final distributions are normalized and smoothed by a convolution with a box filter. We use $\ti=0$ and $\tf = 0.2$ and $\mu=\beta=1$. Assigned boundary conditions (shaded blue area in first and final panels) are given by \eqref{eq:boundary_conditions_initial} with $U_{\iota}(q) = \frac{1}{4}(q-1)^4$ and \eqref{eq:boundary_conditions_final} with $U_{f}(q) = \frac{1}{4}(q^2-1)^2$.}
 \end{figure}

With these results, we have a numerical approximation of the drift which maps an assigned initial probability density into an assigned final density while minimising the Kullback-Liebler divergence on the interval $[\ti,\tf]$. We use this drift to compute the solution of the Fokker-Planck~\eqref{eq:fpk_coupled} $\mathtt{p}_t$ via Algorithm~\ref{alg:fpkgirsanov}, and compare it to the density resulting from the iteration method of \cite{CaHa2022} in Fig.~\ref{fig:overdamped_fpk}.

\section{Fokker-Planck for a time-dependent mechanical underdamped diffusion process}
\label{sec:ud}

We now turn our attention to the underdamped dynamics \cite{ZwaR2001}
\begin{equation}
	\begin{split}
	&\mathrm{d}\bm{\mathscr{q}}_{t}=\frac{\bm{\mathscr{p}}_{t}}{m}\,
	\mathrm{d}t
	\\
	&
	\mathrm{d}\bm{\mathscr{p}}_{t}=-\left(\frac{\bm{\mathscr{p}}_{t}}{\tau}+(\bm{\partial}U_{t})(\bm{\mathscr{q}}_{t})\right)\,\mathrm{d}t
	+\sqrt{ \frac{2\,m}{\tau\,\beta}}\,\mathrm{d}\bm{\mathscr{w}}_{t}		
	 \end{split}
	 \label{ud:sde}
\end{equation}
This is probably the most popular model of an open classical system in contact with a bath at temperature $\beta^{-1}$, driven by a mechanical force and subject to a linear friction force dissipating energy at Stokes rate $\tau$. As in Section~\ref{sec:od}, we assume that the mechanical potential is confining. The corresponding Fokker-Planck equation
\begin{align}
		&
	\left(\partial_{t}+\left \langle\,\frac{\bm{p}}{m}\,,\partial_{\bm{q}} \,\right\rangle-\left \langle\,(\bm{\partial}U_{t})(\bm{q})\,,\partial_{\bm{p}} \,\right\rangle
	-\left \langle\,\bm{\partial}_{\bm{p}}\,,\frac{\bm{p}}{\tau}\,\right\rangle-\frac{ m}{\beta\,\tau}\bm{\partial}_{\bm{p}}^{2}\right)\mathtt{p}_{t}(\bm{x})=0
	\label{ud:FP}
\end{align}
relaxes to a Maxwell-Boltzmann equilibrium for any time independent confining potential
\begin{align}
	\mathtt{p}_{\mathrm{eq}}(\bm{x})=\frac{e^{-\frac{\beta\left\|\bm{p}\right\|^{2}}{2\,m}-\beta\,U(\bm{q})}}{Z}
	\nonumber
\end{align}
with $Z$ a normalizing constant. Our aim is to express the solution of (\ref{ud:FP}) as a suitable Monte Carlo average over an initial datum at time $\ti$.
We proceed analogously to Section~\ref{sec:od} and posit
\begin{align}
	\mathtt{p}_{t}(\bm{x})=e^{-\frac{\beta\left\|\bm{p}\right\|^{2}}{2\,m}}\mathtt{f}_{t}(\bm{x})
	\nonumber
\end{align}
where $\bm{x}=[\bm{q},\bm{p}]^{\top}$. The symplectic component of the drift acts in the same manner on the probability density and the auxiliary function $\mathtt{f}_{t}$. The dissipative component changes sign:
\begin{align}
	\left(\partial_{t}+\left \langle\,\frac{\bm{p}}{m}\,,\partial_{\bm{q}} \,\right\rangle-\left \langle\,(\bm{\partial}U_{t})(\bm{q})-\frac{\bm{p}}{\tau}\,,\bm{\partial}_{\bm{p}}\right\rangle-\frac{ m}{\beta\,\tau}\bm{\partial}_{\bm{p}}^{2}\right)\mathtt{f}_{t}(\bm{x})=-\frac{\beta}{m}\left \langle\,(\bm{\partial} U_{t})(\bm{q})\,,\bm{p} \,\right\rangle\mathtt{f}_{t}(\bm{x})
	\label{ud:aux}
\end{align}
The upshot is that the resulting equation admits the interpretation of a non-homogeneous backward Kolmogorov equation associated to the \textbf{backward} process
\begin{equation}
	\begin{split}
	&\mathrm{d}^{\flat}\bm{\mathscr{q}}_{t}=\frac{\bm{\mathscr{p}}_{t}}{m}
	\,\mathrm{d}t
	\\
	&
	\mathrm{d}^{\flat}\bm{\mathscr{p}}_{t}=\left(\frac{\bm{\mathscr{p}}_{t}}{\tau}-(\bm{\partial}U_{t})(\bm{\mathscr{q}}_{t})\right)\,\mathrm{d}t
	+\sqrt{ \frac{2\,m}{\tau\,\beta}}\,\mathrm{d}\bm{\mathscr{w}}_{t}		
	\end{split}
	\label{ud:bsde}
\end{equation}
Specifically, it is possible to prove that
\begin{proposition}
	The solution of the Fokker-Planck equation (\ref{ud:FP}) can be couched into a conditional expectation
	\begin{align}
		\mathtt{p}_{t}(\bm{x})=e^{-\frac{\beta\left\|\bm{p}\right\|^{2}}{2\,m}}
		\operatorname{E}_{\mathcal{P}^{\flat}}\left(	e^{\frac{\beta\left\|\bm{\mathscr{p}}_{\ti}\right\|^{2}}{2\,m}}\mathtt{p}_{\ti}(\bm{\mathscr{x}}_{\ti})\,e^{-\frac{\beta}{m}\int_{\ti}^{t}\mathrm{d}s\,\left \langle\,(\bm{\partial}U_{s})(\bm{\mathscr{q}}_{s})\,,\bm{\mathscr{p}}_{s} \,\right\rangle}\,\bigg{|}\,\bm{\mathscr{x}}_{t}=\bm{x}\right)
		\label{ud:MC}
	\end{align}
	with respect to the path measure $\mathcal{P}^{\flat}$ generated by (\ref{ud:bsde}).
\end{proposition}
\begin{idea}
	The proof mirrors that of the overdamped case. We define the auxiliary function
	\begin{align}
		\mathtt{g}_{t}(\bm{\mathscr{x}}_{t})=e^{-\frac{\beta}{m}\int_{t}^{\tf}\mathrm{d}s\,\left \langle\,(\bm{\partial}U_{s})(\bm{\mathscr{q}}_{s})\,,\bm{\mathscr{p}}_{s} \,\right\rangle}\mathtt{f}_{t}(\bm{\mathscr{x}}_{t})
		\nonumber
	\end{align}
	for any $\mathrm{f}_t$ solution of (\ref{ud:aux}). Differentiation backward in time along the paths of (\ref{ud:bsde}) yields
	\begin{align}
		\mathrm{d}^{\flat}\mathtt{g}_{t}(\bm{\mathscr{x}}_{t})=
		e^{-\frac{\beta}{m}\int_{t}^{\tf}\mathrm{d}s\,\left \langle\,(\bm{\partial}U_{s})(\bm{\mathscr{q}}_{s})\,,\bm{\mathscr{p}}_{s} \,\right\rangle}\sqrt{ \frac{2\,m}{\tau\,\beta}}\left \langle\,\mathrm{d}^{\flat}\bm{\mathscr{w}}_{t}\,,\bm{\partial}_{\bm{\mathscr{p}}_{t}}\mathtt{g}_{t}(\bm{\mathscr{x}}_{t})\,\right\rangle
		\nonumber
	\end{align}
	We conclude that in any time integral where $\mathrm{g}_{t}$ is integrable, then it is also a martingale with respect to the measure
	$\mathcal{P}^{\flat}$ generated by (\ref{ud:bsde}):
	\begin{align}
		\operatorname{E}_{\mathcal{P}^{\flat}}\left(\mathtt{g}_{\ti}(\bm{\mathscr{x}}_{\ti})\,\Big{|}\,\bm{\mathscr{x}}_{\tf}=\bm{x}\right)
		=
		\operatorname{E}_{\mathcal{P}^{\flat}}\left(\mathtt{g}_{\tf}(\bm{\mathscr{x}}_{\tf})\,\Big{|}\,\bm{\mathscr{x}}_{\tf}=\bm{x}\right)
		=\mathtt{f}_{\tf}(\bm{x})
		\nonumber
	\end{align}
	The chain of identities yields the claim.
\end{idea}

In the case of particular physical interest, when the system probability density at time $\ti$ takes the Maxwell-Boltzmann form
\begin{align}
	\mathtt{p}_{\ti}(\bm{x})=\frac{e^{-\frac{\beta\left\|\bm{p}\right\|^{2}}{2\,m}-\beta\,U_{\iota}(\bm{q})}}{Z}
	\nonumber
\end{align}
then (\ref{ud:MC}) reduces to
\begin{align}
	\mathtt{p}_{t}(\bm{x})
		=e^{-\frac{\beta\left\|\bm{p}\right\|^{2}}{2\,m}}\operatorname{E}_{\mathcal{P}^{\flat}}\left(	\frac{e^{-\beta\,U_{\iota}(\bm{\mathscr{q}}_{\ti})}}{Z}\,e^{-\frac{\beta}{m}\int_{\ti}^{t}\,\mathrm{d}s\,\left \langle\,(\bm{\partial}U_{s})(\bm{\mathscr{q}}_{s})\,,\bm{\mathscr{p}}_{s} \,\right\rangle}\,\bigg{|}\,\bm{\mathscr{x}}_{t}=\bm{x}\right)
	\nonumber
\end{align}

\subsection{Representation with respect to a time reversal invariant measure}

We now turn our attention to the path measure $ \mathcal{I}$ generated by the forward 
process 
\begin{subequations}
	\label{tr:sde}
	\begin{align}
		&\label{tr:sde1}\mathrm{d}\bm{\mathscr{q}}_{t}=\frac{\bm{\mathscr{p}}_{t}}{m}\,
		\mathrm{d}t
		\\
		&\label{tr:sde2}
		\mathrm{d}\bm{\mathscr{p}}_{t}=-
		(\bm{\partial}U_{t})(\bm{\mathscr{q}}_{t})\,\mathrm{d}t
		+\sqrt{ \frac{2\,m}{\tau\,\beta}}\,\mathrm{d}\bm{\mathscr{w}}_{t}	
	\end{align}
\end{subequations}
The drift of the process is divergence-less. This, together with the statistical invariance under time reversal of the Wiener process, imply that we can also interpret the equation
\begin{align}
	\left(\partial_{t}+\left \langle\,\frac{\bm{p}}{m}\,,\partial_{\bm{q}} \,\right\rangle-\left \langle\,(\bm{\partial}U_{t})(\bm{q})\,,\partial_{\bm{p}} \,\right\rangle+\frac{ m}{\beta\,\tau}\bm{\partial}_{\bm{p}}^{2}\right) \mathrm{f}_{t}(\bm{q},\bm{p})=0
	\nonumber
\end{align}
both as a backward Kolmogorov equation or as the Fokker-Planck equation associated to  (\ref{tr:sde}) if we replace forward differentials with \textbf{backward} differentials.
This means that both $\mathcal{P}$ and $\mathcal{P}^{\flat}$ are absolutely continuous with respect to $\mathcal{I}$. 
In particular, we find that
\begin{align}
	\frac{\mathrm{d}\mathcal{P}^{\flat}}{\mathrm{d} \mathcal{I}}=\exp\int_{\ti}^{t}\left(\sqrt{ \frac{\tau\,\beta}{2\,m}}\left \langle\,\mathrm{d}\bm{\mathscr{w}}_{s}\overset{\circ}{\,,}\,\frac{\bm{\mathscr{p}}_{s}}{\tau}\,\right\rangle-\frac{\tau\,\beta}{4\,m}\,\mathrm{d}s\,\left\|\frac{\bm{\mathscr{p}}_{s}}{\tau}\right\|^{2}\right)
	\label{ud:Girsanov}
\end{align}
As in Section~\ref{sec:od}, the symbol $\circ$ denotes that the post-point prescription in the construction of the stochastic integral:
\begin{align}
	\int_{\ti}^{t}\left \langle\,\mathrm{d}\bm{\mathscr{w}}_{s}\overset{\circ}{\,,}\,\frac{\bm{\mathscr{p}}_{s}}{\tau}\,\right\rangle
	=
	\sum_{\ti<\mathfrak{t}_{i} \,\leq\,t} \left \langle\,\bm{\mathscr{w}}_{\mathfrak{t}_{i}}-\bm{\mathscr{w}}_{\mathfrak{t}_{i-1}}\,,\frac{\bm{\mathscr{p}}_{\mathfrak{t}_{i}}}{\tau}\,\right\rangle
	\nonumber
\end{align}
The immediate consequence is the representation of the solution of the Fokker-Planck equation (\ref{ud:FP}) 
\begin{align}
	\mathtt{p}_{t}(\bm{x})=e^{-\frac{\beta\left\|\bm{p}\right\|^{2}}{2 m}}\operatorname{E}_{\mathcal{I}}\left(	\frac{e^{-\beta\,U_{\iota}(\bm{\mathscr{q}}_{\ti})}}{Z}\,e^{-\frac{\beta}{m}\int_{\ti}^{t}\mathrm{d}s\,\left \langle\,(\bm{\partial}U_{s})(\bm{\mathscr{q}}_{s})\,,\bm{\mathscr{p}}_{s} \,\right\rangle}\frac{\mathrm{d}\mathcal{P}^{\flat}}{\mathrm{d} \mathcal{I}}\,\bigg{|}\,\bm{\mathscr{x}}_{t}=\bm{x}\right)
	\label{ud:Imes}
\end{align}
Some further simplifications are possible. In view of the identities
\begin{align}
	\int_{\ti}^{t}\mathrm{d}s\,\left \langle\,(\bm{\partial}U_{s})(\bm{\mathscr{q}}_{s})\,,\bm{\mathscr{p}}_{s} \,\right\rangle
	=
	\int_{\ti}^{t}\left \langle\,\sqrt{ \frac{2\,m}{\tau\,\beta}}\mathrm{d}\bm{\mathscr{w}}_{s}	-\mathrm{d}\bm{\mathscr{p}}_{s}\overset{\circ}{\,,}\,\bm{\mathscr{p}}_{s}\,\right\rangle
	\nonumber
\end{align}
and
\begin{align}
	\mathrm{d}\bm{\mathscr{p}}_{t}^{2}=2\,\left \langle\,\bm{\mathscr{p}}_{t}\overset{\circ}{\,,}\,\mathrm{d}\bm{\mathscr{p}}_{t}\,\right\rangle-\frac{2\,m\,d}{\tau\,\beta}\,\mathrm{d}t
	\nonumber
\end{align}
we can couch (\ref{ud:Imes}) into the form
\begin{align}
	\mathtt{p}_{t}(\bm{x})=e^{\frac{(t-\ti) d}{\tau}}\operatorname{E}_{\mathcal{I}}\left(	\frac{e^{-\beta\,U_{\iota}(\bm{\mathscr{q}}_{\ti})-\frac{\beta\left\|\bm{\mathscr{p}}_{\ti}\right\|^{2}}{2\,m} -\int_{\ti}^{t}\left(\sqrt{ \frac{\beta\,\tau}{2\,m}}\left \langle\,\mathrm{d}\bm{\mathscr{w}}_{s}	\overset{\circ}{\,,}\,\frac{\bm{\mathscr{p}}_{s}}{\tau}\,\right\rangle+\frac{\tau\,\beta}{4\,m}\mathrm{d}s\,\left\|\frac{\bm{\mathscr{p}}_{s}}{\tau}\right\|^{2}\right)}}{Z}\,\Bigg{|}\,\bm{\mathscr{x}}_{t}=\bm{x}\right)
	\nonumber
\end{align}
Finally, we can re-write the stochastic integral  in the pre-point discretization using the identity
\begin{align}
	\int_{\ti}^{t}\sqrt{ \frac{\tau\,\beta}{2\,m}}\left \langle\,\mathrm{d}\bm{\mathscr{w}}_{s}	\overset{\circ}{\,,}\,\frac{\bm{\mathscr{p}}_{s}}{\tau}\,\right\rangle
	=
	\int_{\ti}^{t}\sqrt{ \frac{\tau\,\beta}{2\,m}}\left \langle\,\mathrm{d}\bm{\mathscr{w}}_{s}\,,\frac{\bm{\mathscr{p}}_{s}}{\tau}\,\right\rangle
	+
	\frac{(t-\ti) d}{\tau}
	\nonumber
\end{align}
Thus we arrive at
\begin{align}
	\mathtt{p}_{t}(\bm{x})=\operatorname{E}_{\mathcal{I}}\left(	\frac{1}{Z}\,e^{-\beta\,U_{\iota}(\bm{\mathscr{q}}_{\ti})-\frac{\beta\left\|\bm{\mathscr{p}}_{\ti}\right\|^{2}}{2\,m}-\int_{\ti}^{t}\bigl(\sqrt{ \frac{\beta\,\tau}{2\,m}}\left \langle\,\mathrm{d}\bm{\mathscr{w}}_{s}	\,,\,\frac{\bm{\mathscr{p}}_{s}}{\tau}\,\right\rangle+\frac{\tau\,\beta}{4\,m}\,\mathrm{d}s\,\left\|\frac{\bm{\mathscr{p}}_{s}}{\tau}\right\|^{2}\bigr)}\bigg{|}\bm{\mathscr{x}}_{t}=\bm{x}\right)
	\nonumber
 \end{align}
The advantage of taking averages with respect to the measure $\mathcal{I}$ is that it allows us to use paths of (\ref{tr:sde}) to integrate both the Fokker-Planck equation and a coupled Hamilton-Jacobi-Bellman equation, necessary when constructing numerical solutions of Schr\"odinger bridge problems.

\subsection{Numerical Example}\label{sec:num_fpk_ud}
In this section, we illustrate the numerical integration of a Fokker-Planck equation \eqref{ud:FP} governing the evolution of the joint density of the momentum and position processes following underdamped dynamics~\eqref{ud:FP}. We once again consider the optimal control problem of minimizing the Kullback-Leibler divergence from a free diffusion \eqref{eq:kl_divergence} on the interval $[\ti,\tf]$. This optimal control problem is formulated as the minimization of a Bismut-Pontryagin functional, resulting in the coupled partial differential equations

\begin{subequations}
    \begin{equation}\label{num_ud_fpk}\begin{split}
            \Bigl(\partial_{t}+\bigl\langle\,\frac{\bm{p}}{m}\,,\partial_{\bm{q}} \,\bigr\rangle-\langle\,(\bm{\partial}U_{t})(\bm{q})\,,\partial_{\bm{p}} \,\rangle	- \langle\,\bm{\partial}_{\bm{p}}\,,\frac{\bm{p}}{\tau}\,\rangle-\frac{ m}{\beta\,\tau}\bm{\partial}_{\bm{p}}^{2}\Bigr)\,\mathtt{p}_{t}(\bm{x})&=0
     \end{split}\end{equation}
     \begin{equation}\label{num_ud_hjb}\begin{split}
           \left (\partial_{t}+  \langle\frac{\bm{p}}{m}\,,\,\bm{\partial}_{\bm{p}}\rangle - \langle\bigl(\frac{\bm{p}}{\tau}+(\bm{\partial}U_t)(\bm{q})\bigr)\,,\,\bm{\partial}_{\bm{q}}\rangle+
		\frac{m}{\beta\, \tau}\,\bm{\partial}^2_{\bm{p}}\right )V_{t}(\bm{x})&=- \frac{\beta m}{4\,\tau^2}\|(\bm{\partial}U_t)(\bm{q})\|^2
     \end{split}\end{equation}
\end{subequations}
with the stationarity condition
\begin{align}\label{num_ud_stationarity_condition}
\partial_{\bm{q}} U_t(\bm{q}) = \int_{\mathbb{R}^d}\dd^d \bm{p} \,\frac{\mathrm{p}_t(\bm{p},\bm{q})}{\int_{\mathbb{R}^d}\dd^d \bm{p}\, \mathrm{p}_t(\bm{p},\bm{q})}\,\partial_{\bm{p}} V_t(\bm{p},\bm{q})
\end{align}
We identify \eqref{num_ud_fpk} as the Fokker-Planck and \eqref{num_ud_hjb} as the Hamilton-Jacobi-Bellman equation. Using Girsanov theorem, we find an expression for the intermediate density, the solution of \eqref{num_ud_fpk}, as an expectation

\begin{equation}
    \mathrm{p}_t(\bm{x}) = \operatorname{E}\biggl(P_\iota(\bm{\mathscr{p}}_{\ti}\,,\, \bm{\mathscr{q}}_{\ti})e^{-\frac{\tau\beta}{4\,m}\int_{\ti}^{t}\dd s \|\bm{\mathscr{p}}_s\|^2- \sqrt{\frac{\tau\beta}{2\,m}}\int_{\ti}^{t}\langle\dd \bm{\mathscr{w}}_s\,,\,\bm{\mathscr{p}}_s\rangle}\,\bigg{|}\,\bm{\mathscr{x}}_t = \bm{x}\biggr)
\end{equation}
taken over trajectories of the process~\eqref{tr:sde}.

 We consider the case where the boundary conditions are assigned on the joint distribution at initial and final time
\begin{subequations}
    \begin{align}
    \label{UD:initial_boundarycond}
    P_{\iota}(\bm{p},\bm{q}) &= \dfrac{1}{Z_\iota}\exp\Biggl(-\dfrac{\beta\|\bm{p}\|^2}{2} - \beta\,U_{\iota}(\bm{q})\Biggr) \\\label{UD:final_boundarycond}
    P_{f}(\bm{p},\bm{q}) &= \dfrac{1}{Z_f}\exp\Biggl(-\dfrac{\beta\|\bm{p}\|^2}{2} - \beta \,U_f(\bm{q})\Biggr)\end{align}
\end{subequations}

with $Z_\iota,Z_f$ normalizing constants. 

\begin{figure}
     \centering
    \includegraphics[width=\linewidth]{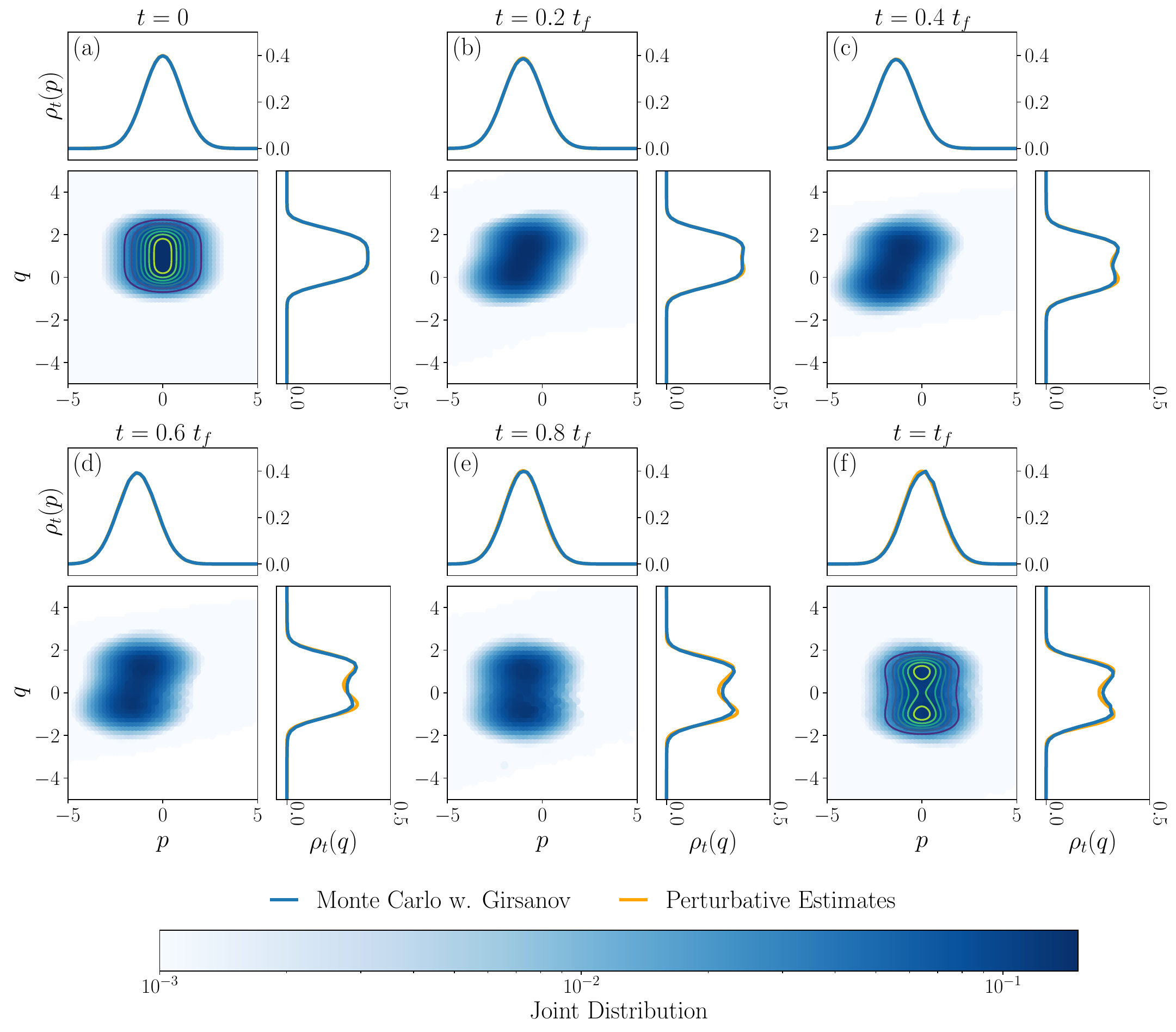}
.     \caption{\label{fig:fpk_underdamped} Solution of a Fokker-Planck equation driven by a non-linear mechanical underdamped diffusion computed by Monte Carlo integration. Panels (a)-(f) show the following: Center: the joint distribution for the momentum and position; Top: the marginal distribution of the momentum; Left: marginal distribution of the position. The optimal protocol $U$ used in the integration and the reference solutions for the marginal densities (orange) are estimated from a perturbative expansion around the overdamped limit \cite{SaBaMG2024}. \newline For the integration, we evolve $M = 10\,000$ trajectories from a set of $2601$ equally spaced points from the interval $[-5,5]\times [-5,5]$. We use a time step size $h=0.025$ and integrate over trajectories of \eqref{ud:FP} using an Euler-Maruyama discretization. We use $\ti=0$, $\tf = 5$, $\beta=25$ and $\tau=m=1$. 
The assigned initial condition is given by \eqref{UD:initial_boundarycond} with $U_{\iota}(q) = \frac{1}{4}(q-1)^4$ and final condition \eqref{UD:final_boundarycond} with $U_{f}(q) = \frac{1}{4}(q^2-1)^2$.}
 \end{figure}

The numerical computation is summarized in Algorithm~\ref{alg:UD_fpkgirsanov}. For the optimal control potential and benchmark solution, we use numerical predictions from \cite{SaBaMG2024}. There, predictions are made for the optimal protocol in the underdamped dynamics using a multiscale perturbative expansion around the overdamped problem; for more detail, see Section 8.2 of \cite{SaBaMG2024}. The prediction for the optimal control protocol is used as the drift in the integration of the Fokker-Planck~\eqref{num_ud_fpk}, with the results shown in Fig.~\ref{fig:fpk_underdamped}.

\begin{algorithm}
\caption{Integrating a Fokker-Planck for an underdamped diffusion process using Girsanov theorem}
\label{alg:UD_fpkgirsanov}
\begin{algorithmic} 
\STATE Initialize $\bm{\mathscr{q}}_{t_n} = \bm{q}\in \mathbb{R}^d$
\STATE Initialize $g = 0$
\STATE Initialize $\partial U_{t_n}$ for $n\in\{0,1,\ldots,N\}$
\STATE Initialize $\delta_{t} = |t_{n-1}-t_n|$
\FOR{ $i$ in $0,\ldots, n-1$}
\STATE Sample Brownian noise: $\epsilon \sim \mathcal{N}(0,1)$
\STATE Evolve one backward step of \eqref{ud:Imes}: $\begin{cases}
    
\bm{\mathscr{q}}_{t_{n-1-i}} = \bm{\mathscr{q}}_{t_{n-i}} - \dfrac{1}{m}\bm{\mathscr{p}}_{t_{i}}\,\delta_t \\[5pt]
\bm{\mathscr{p}}_{t_{n-1-i}} = \bm{\mathscr{p}}_{t_{n-i}} + (\partial U_{t_{n-i}})(\bm{\mathscr{q}}_{t_{n-i}})\,\delta_t - \sqrt{\delta_{t}\frac{2m\,}{\tau\,\beta}}\, \epsilon\end{cases}$ 
\STATE Add Girsanov weight: $g = g+ \delta_t \frac{\tau\beta}{4\,m}   \,\bm{p}_{t_{n-1-i}}^2 +\sqrt{\delta_t\frac{\tau\beta}{2\,m}\,}\,\epsilon\, \bm{p}_{t_{n-1-i}}$
\ENDFOR
\STATE Return $\mathtt{p}_{t_n}(\bm{q}) = P_{\iota}(\bm{\mathscr{p}}_{t_\iota},\bm{\mathscr{q}}_{t_\iota})\,e^{-g}$
\end{algorithmic}
\end{algorithm}

\section{Bismut-Elworthy-Li Monte Carlo representation of gradients}\label{sec:od_bel}

Numerical integration of Schr\"odinger bridge type problems, in the overdamped \cite{AuGaMeMoMG2012,LeoC2014,CaHa2022,ChGePa2021} and underdamped \cite{ PMG2014,MGSc2014,SaBaMG2024,SaBaMG2024_2} cases require the solution of a Hamilton-Jacobi-Bellman (also known as a dynamic programming) equation, specifying the optimal control potential.
In the simplest overdamped setup, the mechanical force is given by (\ref{stationarity_condition}).
The function  
\begin{align}
	V_{t}(\bm{q})\colon[\ti,\tf]\,\times\,\mathbb{R}^{d}\mapsto\mathbb{R}
	\nonumber
\end{align}
solves a Burgers type equation. More generally, optimization problems often require computing gradients of scalar functions satisfying a non-homogeneous backward Kolmogorov equation in $[\ti\,,\tf]$ of the form
\begin{subequations}
	\label{BEL:dp}
	\begin{align}
		\label{BEL:dp1}
		\left (\partial_{t}+ \left \langle\,\bm{b}_{t}(\bm{x})\,,\bm{\partial}_{\bm{x}}\,\right\rangle+
		\left \langle\,\mathsf{A}_{t}(\bm{x})\mathsf{A}_{t}^{\top}(\bm{x}) \,, \bm{\partial}_{\bm{x}}\otimes\bm{\partial}_{\bm{x}}\,\right\rangle\right )V_{t}(\bm{x})&=\,-\,F_{t}(x)
		\\
		\label{BEL:dp2}
		V_{\tf}(\bm{x})&=\varphi(\bm{x})
	\end{align}
\end{subequations}
The left hand side of (\ref{BEL:dp1}) is the mean forward derivative (see Chapter~11 of \cite{NelE2001}) of $V_{t}$ along the paths of the $n$-dimensional system of It\^o stochastic differential equations
\begin{align}
	\mathrm{d}\bm{\mathscr{x}}_{t}=\bm{b}_{t}(\bm{\mathscr{x}}_{t})\,\mathrm{d}t+\mathsf{A}_{t}(\bm{\mathscr{x}}_{t})\,\mathrm{d}\bm{w}_{t}
	\label{BEL:sde}
\end{align}
In (\ref{BEL:dp1}) and (\ref{BEL:sde}), we consider drift $\bm{b}$ and volatility fields $\mathsf{A}$ of more general form than in
Sections~\ref{sec:od} and \ref{sec:ud}. This choice means that the following discussion is applicable to both overdamped and underdamped cases, as well as to more general situations, 
including non-linear problems \cite{EHuJeKr2021}. In non-linear problems, the expression of the solution of (\ref{BEL:dp}) and its gradient are iteratively computed in sequences of infinitesimal time horizons $[\ti,\tf]$ to construct the solution of partial differential equations in which $\bm{b}$, $\mathsf{A}$ and $F$ depend upon the unknown field $V$.

It is well known that Dynkin's formula (see e.g. Chapter~6 of \cite{KleF2005}) yields a Monte Carlo representation of the solution of (\ref{BEL:dp})
\begin{align}
	V_{t}(\bm{x})=\operatorname{E}_{\mathcal{P}}\left( \varphi(\bm{\mathscr{x}}_{\tf})+\int_{t}^{\tf}\mathrm{d}s\,F_{s}(\bm{\mathscr{x}}_{s})
	\,\Big{|}\,\bm{\mathscr{x}}_{t}=\bm{x}\right)
	\label{Dynkin}
\end{align}
Our goal is to find an analogous expression for the gradient of $V_{t}$. The Bismut-Elworthy-Li formula \cite{BisJ1984,ElLi1994,ElLi1994b} accomplishes this task. 
\begin{remark}
	In what follows, to neaten mathematical formulae we adopt the push-forward notation for the Jacobian matrix of a vector field.
	We refer to Section $\mathscr{O}.j$ pag. $\mathrm{xlii}$ of \cite{FraT2012} for a geometrical justification of the notation.
 
 For any $\bm{v}\colon \mathbb{R}^{d}\mapsto\mathbb{R}^{d}$, we write
	\begin{align}
		\left \langle\,\bm{e}_{i}\,, \bm{v}_{*} \bm{e}_{j}\,\right\rangle:= \partial_{x^{j}} v_{t}^{i}(\bm{x})
		\label{BEL:pf}
	\end{align}
	where $ \bm{e}_{i}$  and $ \bm{e}_{j}$ are respectively the $i$-th and $j$-th elements of the canonical basis of $\mathbb{R}^{d}$.
  Under our regularity assumptions, we regard the solution of (\ref{BEL:sde}) satisfying the condition 
  \begin{align}
&  	\bm{\mathscr{x}}_{s}=\bm{x} && s\,\leq\,t
  	\nonumber
  \end{align}
  as the image of the stochastic flow $\bm{X}\colon \mathbb{R}\times \mathbb{R}\times \mathbb{R}^{d}\mapsto \mathbb{R}^{d}$ \cite{KunH1986} such that
  \begin{align}
  	\bm{\mathscr{x}}_{t}=\bm{X}_{t,s}(\bm{\mathscr{x}}_{s})
  	\nonumber
  \end{align}
  and omit reference to the initial data on the left hand side, when no ambiguity arises. 
  According to (\ref{BEL:pf}), we denote the cocycle obtained by differentiating the flow $\bm{X}_{t,s}$ with respect to its argument as $\bm{\mathscr{x}}_{*t,s}$ implying that
	\begin{align}
		\left \langle\,\bm{e}_{i}\,,\bm{\mathscr{x}}_{*t,s}\bm{e}_{j}\,\right\rangle
		:=\partial_{\mathscr{x}^{j}}X_{t,s}^{i}(\bm{x})
		\nonumber
	\end{align}
 By definition  $\bm{\mathscr{x}}_{*t,s}$ enjoys the cocycle property \cite{ArnL1998}, meaning that 
 \begin{align}
& \bm{\mathscr{x}}_{*t,s}(\bm{\mathscr{x}}_{s})=\bm{\mathscr{x}}_{*t,u}(\bm{\mathscr{x}}_{u})
\bm{\mathscr{x}}_{*u,s}(\bm{\mathscr{x}}_{s}) &\forall\, s\le u \le t
	\nonumber
	\end{align} 
	\begin{center}
		* *
	\end{center}
\end{remark}
Here, we present a heuristic, physics style derivation of the formula based on Malliavin's stochastic variational calculus \cite{NuaD1995} which draws from the mathematically more rigorous exposition in \cite{ZhaY2010}, and is close to the original treatment in \cite{BisJ1984}. To this end, we observe that if $\bm{e}_{i}$
is the $i$-th element of the canonical basis of $\mathbb{R}^{n}$
\begin{align}
	\left \langle\bm{e}_{i} \,,(\bm{\partial}	V_{t})(\bm{x})\right\rangle=\operatorname{E} \left( \left  \langle\bm{\mathscr{x}}_{*\tf,t}\bm{e}_{i}\,,(\bm{\partial} \varphi)(\bm{\mathscr{x}}_{\tf})\right \rangle
	+\int_{t}^{\tf}\mathrm{d}s\,
	\left  \langle\bm{\mathscr{x}}_{*\tf,t}\bm{e}_{i}\,,(\bm{\partial} F_{s})(\bm{\mathscr{x}}_{s})\right \rangle
 \Big{|}\bm{\mathscr{x}}_{t}=\bm{x}\right )
	\label{BEL:formula}
\end{align}
where $ \bm{\mathscr{x}}_{*\tf,t}$ denotes the matrix valued process obtained by varying (\ref{BEL:sde}) with respect to its initial
datum. In other words, if we suppose
\begin{align}
&	\bm{\mathscr{x}}_{s}=\bm{x}&& s\,\leq\,t
	\nonumber
\end{align}
then
\begin{align}
&	\mathrm{d}\bm{\mathscr{x}}_{*t,s}=\Big{(}\bm{b}_{*t}(\bm{\mathscr{x}}_{t})\,\mathrm{d}t+\mathsf{A}_{*t}(\bm{\mathscr{x}}_{t})\,\mathrm{d}\bm{w}_{t}\Big{)}\,\bm{\mathscr{x}}_{*t,s}
\nonumber\\
& \bm{\mathscr{x}}_{*s,s}=\operatorname{1}_{n}
	\nonumber
\end{align}
The identity (\ref{BEL:formula}) allows us to derive the Bismut-Elworthy-Li formula from Malliavin's integration by parts formula.

\subsection{Integration by parts formula}

Let us consider the equation
\begin{align}
	\begin{split}
&		\mathrm{d}\bm{\mathscr{x}}_{t}^{(\varepsilon)}=\left(\bm{b}_{t}(\bm{\mathscr{x}}_{t}^{(\varepsilon)})+\varepsilon\,\bm{\mathscr{h}}_{t}\right)\,\mathrm{d}t+\mathsf{A}_{t}(\bm{\mathscr{x}}_{t}^{(\varepsilon)})\,\mathrm{d}\bm{w}_{t}
	\\
&				\bm{\mathscr{x}}_{s}^{(\varepsilon)}=\bm{x}
		\end{split}
\label{M:epseq}
\end{align}
We assume $\bm{\mathscr{h}}_{t}$ to be a differentiable process, although rigorous constructions of integration by parts formula see e.g. \cite{NuaD1995} weaken this assumption to processes of bounded variation (see Chapter~1 of \cite{KleF2005}). Differentiating at $ \varepsilon=0$ yields the variational equation
\begin{align}
	\begin{split}
		&		\mathrm{d}\bm{\mathscr{x}}_{t}^{\prime}=\Big{(}\bm{b}_{*t}(\bm{\mathscr{x}}_{t})\,\mathrm{d}t+\mathsf{A}_{*t}(\bm{\mathscr{x}}_{t})\,\mathrm{d}\bm{w}_{t}\Big{)}\,\bm{\mathscr{x}}_{t}^{\prime}
		+
		\bm{\mathscr{h}}_{t}\,\mathrm{d}t
		\\
		&\bm{\mathscr{x}}_{s}^{\prime}=0
	\end{split}
	\label{M:vareq}
\end{align}
We can always write the solution of this latter equation in terms of the push-forward of the flow of \eqref{BEL:sde}
\begin{align}
	\bm{\mathscr{x}}_{t}^{\prime}=\int_{s}^{t}\mathrm{d}u\,\bm{\mathscr{x}}_{*t,u}\,\bm{\mathscr{h}}_{u}
	\label{M:xprime}
\end{align}
Therefore for sufficiently small $\varepsilon$ 
\begin{align}
	\bm{\mathscr{x}}_{t}^{(\varepsilon)}=\bm{\mathscr{x}}_{t}+\varepsilon\,\bm{\mathscr{x}}_{t}^{\prime}+\mathrm{h.o.t}
	\label{M:function}
\end{align}
allows us to regard the solution of \eqref{M:vareq} as a functional of the solution of \eqref{BEL:sde} (h.o.t. stands for higher order terms).
The conclusion is that we can compute the expectation value of any integrable function $g$ of a solution of (\ref{M:vareq}) by expressing it 
as a function of the solution of (\ref{BEL:sde}) via (\ref{M:function}) and then averaging with respect to the measure $\mathcal{P}$ generated by (\ref{BEL:sde})
\begin{align}
\operatorname{E}_{\mathcal{P}^{\varepsilon}}\left( g(\bm{\mathscr{x}}_{t})\,\big{|}\,\bm{\mathscr{x}}_{s}=\bm{x}\right)=\operatorname{E}_{\mathcal{P}}\left( g(\bm{\mathscr{x}}_{t}^{(\varepsilon)})\,\big{|}\,\bm{\mathscr{x}}_{s}=\bm{x}\right)
	\label{M:left}
\end{align}
A second connection comes from Girsanov's change of measure formula. Namely, if $\mathcal{P}^{\varepsilon}$ is the path measure generated by (\ref{M:epseq}), then for any test function $g$, we get the identity
\begin{align}
&	\operatorname{E}_{\mathcal{P}^{\varepsilon}}\left(g(\bm{\mathscr{x}}_{t})\,\Big{|}\,\bm{\mathscr{x}}_{s}=\bm{x}\right) =\operatorname{E}_{\mathcal{P}}
	\left(g(\bm{\mathscr{x}}_{t})\frac{\mathrm{d}\mathcal{P}^{\varepsilon}}{\mathrm{d} \mathcal{P}}\,\Big{|}\,\bm{\mathscr{x}}_{s}=\bm{x}\right)
	\nonumber\\
&=\operatorname{E}_{\mathcal{P}}
\biggl(g(\bm{\mathscr{x}}_{t})
\,\exp\int_{s}^{t}\biggl(\left \langle\,\mathrm{d}\bm{\mathscr{w}}_{u}\,,\varepsilon\,\mathsf{A}_{u}^{-1}(\bm{\mathscr{x}}_{u})\bm{\mathscr{h}}_{u}\,\right\rangle-\mathrm{d}u\,\frac{\left\|\varepsilon\,\mathsf{A}_{u}^{-1}(\bm{\mathscr{x}}_{u})\,\bm{\mathscr{h}}_{u}\right\|^{2}}{2}
\biggr)\,\Big{|}\,\bm{\mathscr{x}}_{s}=\bm{x}
\biggr)
	\label{M:right}
\end{align}
If $g$ is also sufficiently regular, upon differentiating (\ref{M:left}) and (\ref{M:right}) at $\varepsilon=0$, we arrive at Malliavin's integration by parts formula
\begin{align}
	\operatorname{E}_{\mathcal{P}}\left(\left \langle\,\bm{\mathscr{x}}_{t}^{\prime}\,, (\bm{\partial}g)(\bm{\mathscr{x}}_{t})\,\right\rangle\,\Big{|}\,\bm{\mathscr{x}}_{s}=\bm{x}\right)
	=\operatorname{E}_{\mathcal{P}}\left(g(\bm{\mathscr{x}}_{t})\int_{s}^{t}\left \langle\,\mathrm{d}\bm{\mathscr{w}}_{u}\,,\mathsf{A}_{u}^{-1}(\bm{\mathscr{x}}_{u})\bm{\mathscr{h}}_{u}\,\right\rangle\,\Big{|}\,\bm{\mathscr{x}}_{s}=\bm{x}\right)
	\nonumber
\end{align}

\subsection{Application to non-degenerate diffusion}

We set
\begin{align}
	\bm{\mathscr{h}}_{u}=\bm{\mathscr{x}}_{*u,s}\bm{e}_{i}
\label{M:nd}
\end{align}
where $ \bm{e}_{i}$ is the $i$-th element of the canonical basis of $\mathbb{R}^{n}$. This is legitimate because, under standard regularity assumptions, $\bm{\mathscr{x}}_{*u,s}$ is a process of finite variation. Upon inserting into (\ref{M:xprime}) we get
\begin{align}
	\bm{\mathscr{x}}_{t,s}^{\prime}=(t-s)\,\bm{\mathscr{x}}_{*t,s}\,\bm{e}_{i}
	\nonumber
\end{align}
The integration by parts formula becomes
\begin{align}
	\operatorname{E}_{\mathcal{P}}\left(\left \langle\,\bm{\mathscr{x}}_{*t,s}\bm{e}_{i}\,, (\bm{\partial}g)(\bm{\mathscr{x}}_{t})\,\right\rangle\,\big{|}\, \bm{\mathscr{x}}_{s}=\bm{x}\right)
	=\operatorname{E}_{\mathcal{P}}\left(\frac{g(\bm{\mathscr{x}}_{t})}{t-s}\int_{s}^{t}\left \langle\,\mathrm{d}\bm{\mathscr{w}}_{u}\,,\mathsf{A}_{u}^{-1}(\bm{\mathscr{x}}_{u})\bm{\mathscr{x}}_{*u,s}\bm{e}_{i}\,\right\rangle\,\Big{|}\, \bm{\mathscr{x}}_{s}=\bm{x}\right)
	\nonumber
\end{align}
The identity holds for arbitrary $t\,\geq\,s$. Hence, we can apply it to (\ref{BEL:formula}) in order to derive the expression of the gradient of the solution of  (\ref{BEL:dp}) 
\begin{align}
    \label{BEL:BEL}
        \left \langle\bm{e}_{i} \,,(\bm{\partial}	V_{t})(\bm{x})\right\rangle&=
	\operatorname{E} \left( \frac{\varphi(\bm{\mathscr{x}}_{\tf})}{\tf-t} \int_{t}^{\tf}\left \langle\,\mathrm{d}\bm{\mathscr{w}}_{u}\,,\mathsf{A}_{u}^{-1}(\bm{\mathscr{x}}_{u})\bm{\mathscr{x}}_{*u,t}\bm{e}_{i}\,\right\rangle
	\,\Big{|}\,\bm{\mathscr{x}}_{t}=\bm{x}\right )
	\nonumber\\
&+	\operatorname{E} \left( 
\int_{t}^{\tf}\mathrm{d}s\,\frac{F_{s}(\bm{\mathscr{x}}_{s})}{s-t}\int_{t}^{s}\left \langle\,\mathrm{d}\bm{\mathscr{w}}_{u}\,,\mathsf{A}_{u}^{-1}(\bm{\mathscr{x}}_{u})\bm{\mathscr{x}}_{*u,t}\bm{e}_{i}\,\right\rangle\,\Big{|}\,\bm{\mathscr{x}}_{t}=\bm{x}\right )
\end{align}
provided the volatility field $\mathsf{A}$ is always non-singular.

\subsubsection{Application to the transition probability density}

It is worth noticing the following consequence of (\ref{BEL:BEL}) when $F_{s}$ vanishes. In such a case, (\ref{BEL:BEL})
reduces to
\begin{align}
	\bm{\partial}_{\bm{x}}\int_{\mathbb{R}^{d}}\mathrm{d}^{d}\bm{y} \,\varphi(\bm{y})\,\operatorname{p}_{\tf,t}(\bm{y}\,|\,\bm{x})
	=\frac{1}{\tf-t}\,\operatorname{E}_{\mathcal{P}}\left(\varphi(\bm{\mathscr{x}}_{\tf})
	\int_{t}^{\tf}\left \langle\,\mathrm{d}\bm{w}_{u}\,,\mathsf{A}_{u}^{-1}(\bm{\mathscr{x}}_{u}))\bm{\mathscr{x}}_{*\,u,t}\right\rangle\,\Big{|}\,\bm{\mathscr{x}}_{t}=\bm{x}\right)
	\nonumber
\end{align}
As the identity must hold true for any $\varphi$,we can also write 
\begin{align}
	 	\bm{\partial}_{\bm{x}}\operatorname{p}_{\tf,t}(\bm{y}\,|\,\bm{x})
	=\frac{1}{\tf-t}\,\operatorname{E}_{\mathcal{P}}\left(\delta^{(d)}(\bm{y}-\bm{\mathscr{x}}_{\tf})
	\int_{t}^{\tf}\left \langle\,\mathrm{d}\bm{w}_{u}\,,\mathsf{A}_{u}^{-1}(\bm{\mathscr{x}}_{u})\bm{\mathscr{x}}_{*\,u,t}\right\rangle
	\,\Big{|}\,\bm{\mathscr{x}}_{t}=\bm{x}
	\right)
	\label{BEL:delta}
\end{align}
 A result by Molchanov, Section~5 of \cite{MolS1975}, allows us to express (\ref{BEL:delta}) in terms of an expectation value with respect to a reciprocal process, see e.g. \cite{KreA1997}. Namely, given a Markov process in $[\ti,\tf]$, we can use it to construct a reciprocal process, i.e. a process conditioned at both ends of the time horizon from the relations
\begin{align}
&	\mathtt{p}_{t,\tf,\ti}(\bm{x}\mid\bm{z},\bm{y})=\frac{\mathtt{p}_{\tf,t}(\bm{z}\mid\bm{x})\,\mathtt{p}_{t,\ti}(\bm{x}\mid\bm{y})}{\mathtt{p}_{\tf,\ti}(\bm{z}\mid\bm{y})}&& \ti \,\leq\,t\,\leq\,\tf
\nonumber\\
&\mathtt{p}_{t_{\mathfrak{2}},t_{\mathfrak{1}}\tf,\ti}(\bm{x}_{\mathfrak{2}},\bm{x}_{\mathfrak{1}}\mid\bm{z},\bm{y})=\frac{\mathtt{p}_{\tf,t_{\mathfrak{2}}}(\bm{z}\mid\bm{x}_{\mathfrak{2}})\,\mathtt{p}_{t_{\mathfrak{2}},t_{\mathfrak{1}}}(\bm{x}_{\mathfrak{2}}\mid\bm{x}_{\mathfrak{1}})\,\mathtt{p}_{t_{\mathfrak{1}},\ti}(\bm{x}_{\mathfrak{1}}\mid\bm{y})}{\mathtt{p}_{\tf,\ti}(\bm{z}\mid\bm{y})}&& \ti \,\leq\,t_{\mathfrak{1}}\,\leq\,t_{\mathfrak{2}}\,\leq\,\tf
\nonumber\\
& \mathrm{etc}
	\label{BEL:reciprocal}
\end{align}
Upon contrasting (\ref{BEL:reciprocal}) with (\ref{BEL:delta}), we thus arrive at Bismut's formula (pag. 78 of \cite{BisJ1984}) for the gradient of the transition probability density
\begin{align}
	 &	\frac{\bm{\partial}_{\bm{x}}\operatorname{p}_{\tf,t}(\bm{y\mid\bm{x})}}{\operatorname{p}_{\tf,t}(\bm{y}\mid\bm{x})}
	=\operatorname{E}_{\mathcal{P}}\left(	\int_{t}^{\tf}\left \langle\,\mathrm{d}\bm{w}_{u}\,,\mathsf{A}_{u}^{-1}(\bm{\mathscr{x}}_{u})\bm{\mathscr{x}}_{*\,u,t}\right\rangle\,\Big{|}\,\bm{\mathscr{x}}_{\tf}=\bm{y},\bm{\mathscr{x}}_{t}=\bm{x}\right)
	\label{BEL:Bismut}
\end{align}
The underscript $\mathcal{P}$ here means that we construct all finite dimensional approximations to the reciprocal process from the transition probability density of (\ref{BEL:sde}) according to (\ref{BEL:delta}).

Unfortunately, (\ref{BEL:Bismut}) does not directly provide a Monte Carlo representation of the score function because the derivative acts on the variable expressing the condition. It is, however, possible to use ideas similar to these and from the previous sections to obtain a Monte Carlo representation of the score function. 

\subsection{Analytical Example}\label{sec:an_hjb_od}

It is worth illustrating the use of Dynkin's and Bismut-Elworthy-Li formulas in a case where all calculations can be performed explicitly.
To this end let us consider 
\begin{align}
\begin{split}
&	\mathrm{d}\mathscr{q}_{t}=\frac{\mathscr{p}_{t}}{m}+\sqrt{\eta\frac{2\,\tau}{m\,\beta}}\,\mathrm{d}w_{t}^{(\mathfrak{1})}
\\
&\mathrm{d}\mathscr{p}_{t}=-\frac{\mathscr{p}_{t}}{\tau}+\sqrt{\frac{2\,m}{\beta\,\tau}}\,\mathrm{d}w_{t}^{(\mathfrak{2})}
	\end{split}
 \label{exa:sde}
\end{align}
whose solution is simply
\begin{align}
&\mathscr{q}_{t}=\mathscr{q}_{\ti}+\sqrt{\eta\frac{2\,\tau}{m\,\beta}}\,w_{t}^{(\mathfrak{1})}+\int_{\ti}^{t}\mathrm{d}s\,\frac{\mathscr{p}_{s}}{m}
\nonumber\\
&	\mathscr{p}_{t}=\mathscr{p}_{\ti}e^{-\frac{t}{\tau}} + \sqrt{\frac{2\,m}{\beta\,\tau}}\int_{\ti}^{t}\mathrm{d}w_{s}^{(\mathfrak{2})}\,e^{-\frac{t-s}{\tau}}
	\nonumber
\end{align}
Let us consider the partial differential equation
\begin{align}
\begin{split}
	&	\partial_{t}V(q,p)+\frac{p}{m}\,\partial_{q}V(q,p)-\frac{p}{\tau}\,\partial_{p}V(q,p)+\frac{\eta\,\tau}{m\,\beta}\,\partial_{q}^{2}V(q,p)+\frac{m}{\beta\,\tau}\,\partial_{p}^{2}V(q,p)=0
	\\
	& V_{\tf}(q,p)=p
	\end{split}
 \label{exa:pde}
\end{align}
It is straightforward to verify that at any time $t\le \tf$
\begin{align}
	V_{t}(p,q)=p\, e^{-\frac{\tf-t}{\tau}}
	\nonumber
\end{align}
Upon applying Dynkin's formula (\ref{Dynkin}), we verify that
\begin{align}
	V_{t}(p,q)=\operatorname{E}\left (\mathscr{p}_{\tf}\,\big{|}\,\mathscr{q}_{t}=q,\mathscr{p}_{t}=p\right )
	\nonumber
\end{align}
Next, we wish to apply Bismut-Elworthy-Li to recover
\begin{align}
	\partial_{p}V_{t}(p,q)=e^{-\frac{\tf-t}{\tau}}
	\label{exa:BEL}
\end{align}
To this end, we determine the co-cycle solution of the linearized dynamics. The co-cycle equation is
\begin{align}
&	\dot{\mathscr{x}}_{*t,s}=\begin{bmatrix}
	0	& \frac{1}{m} \\
	0  &-\frac{1}{\tau}
	\end{bmatrix}\mathscr{x}_{*t,s}
&&	\mathscr{x}_{*s,s}=\begin{bmatrix}
	1	& 0 \\
	0  & 1
	\end{bmatrix}
	\nonumber
\end{align}
from where we get the unique solution
\begin{align}
	\mathscr{x}_{*t,s}=\begin{bmatrix}
	1 & \frac{1}{m}\Bigl(1-e^{-\frac{t-s}{\tau}}	\Bigr)\\  0 & e^{-\frac{t-s}{\tau}}
	\end{bmatrix}
	\nonumber
\end{align}
To evaluate Bismut-Elworthy-Li formula, we also need the inverse of the volatility matrix, which is
\begin{align}
	\operatorname{A}^{-1}=\begin{bmatrix}
	\sqrt{\frac{\beta\,m}{2\,\eta \,\tau}}	& 0\\  0 & \sqrt{\frac{\beta\,\tau}{2\,m}}
	\end{bmatrix}
 \nonumber
\end{align}
We thus obtain
\begin{align}
\partial_{p}V_{t}(p,q)&=
\operatorname{E}\Biggl(\frac{\mathscr{p}_{t}\,e^{-\frac{\tf-t}{\tau}}+\sqrt{\frac{2\,m}{\beta\,\tau}}\int_{t}^{\tf}\,\mathrm{d}w_{s}^{(\mathfrak{2})}\,e^{-\frac{\tf-s}{\tau}}}{\tf -t} \times 
\nonumber\\
&\biggl(
	\int_{t}^{\tf}
	\mathrm{d}w_{u}^{(\mathfrak{1})}\sqrt{\frac{\beta\,m}{2\,\tau}}	\,\frac{1-e^{-\frac{u-t}{\tau}}}{m}
	+\sqrt{\frac{\beta\,\tau}{2\,m}}\,\mathrm{d}w_{u}^{(\mathfrak{2})} e^{-\frac{u-t}{\tau}}
	\biggr)\,
	\bigg{|}\,\mathscr{q}_{t}=q,\,\mathscr{p}_{t}=p\Biggr)
 \nonumber
\end{align}
Using standard properties of stochastic integrals \cite{KleF2005}, we recover the expected result
\begin{align}
\partial_{p}V_{t}(p,q)=\int_{t}^{\tf}\mathrm{d}s\,\frac{e^{-\frac{1}{\tau}(\tf-s)}e^{-\frac{1}{\tau}(s-t)}}{\tf -t}=e^{-\frac{\tf-t}{\tau}}
\nonumber
\end{align}

\subsection{Numerical Example}\label{sec:num_hjb_od}
In this section, we apply the Bismut-Elworthy- Li to compute the gradient of the value function in the optimal control problem minimizing the Kullback-Leibler divergence \eqref{eq:kl_divergence} in the overdamped dynamics: the gradient of the solution to \eqref{eq:hjb_coupled}. This is calculated as a numerical average over sampled trajectories of \eqref{od:FP}. We use the same approximation of the optimal control potential $U$ as in the Fokker-Planck example of Section~\ref{sec:num_fpk_od}. We find 
\begin{equation}
 \bm{\mathscr{q}}_{*s,t} = \begin{cases}
          1 &\text{ for } s=t
 \\[5pt]
e^{-\mu \int_t^s \dd r \,\partial^2 U_r(\bm{\mathscr{q}}_r)}
& s\,>\,t
 \end{cases}
\end{equation}
Hence \eqref{BEL:BEL} becomes 
\begin{equation}
    \begin{split}
	 &\,(\partial_{\bm{q}}	V_{t})(\bm{q}) =
	\sqrt{\frac{\beta}{2\,\mu}}\operatorname{E}\Biggl( \frac{\varphi(\bm{\mathscr{q}}_{\tf})}{\tf - t}\int_t^{\tf}\Bigl\langle\dd \bm{\mathscr{w}}_s\, , \,\exp\Bigl(-\mu \int_t^{s}\dd r\,\partial^2 U_r(\bm{\mathscr{q}}_r)\Bigr)\Bigr\rangle \\
	 &+ \frac{\beta\,\mu}{4}
\int_{t}^{\tf}\mathrm{d}s\,\dfrac{{\|}\partial U_{s}(\bm{\mathscr{q}}_{s}){\|}^2}{s-t}\,\int_t^{s}\Bigl\langle\dd \bm{\mathscr{w}}_r\,, \,\exp\Bigl(-\mu \int_t^{r}\dd v\,\partial^2 U_v(\bm{\mathscr{q}}_v)\Bigr)\Bigr\rangle \,\bigg{|}\,\bm{\mathscr{q}}_{t}=\bm{q}\Biggr)
	\label{BEL:BEL_solution}
    \end{split}
\end{equation}
We repeatedly sample trajectories of the process \eqref{od:FP} using the Euler-Maruyama discretization scheme and compute the integrals as running costs over each trajectory, finally taking a numerical expectation. The calculation is summarized in Algorithm~\ref{alg_dvbel_over}.

From the physics point of view, note that we can conceptualize the motility constant with the ratio
\begin{align} \mu= \frac{\tau}{m} \end{align}
for consistency with the underdamped equations.

\begin{algorithm}
\caption{Monte Carlo integration for gradient of the value function}
\label{alg_dvbel_over}
\begin{algorithmic} 
\STATE Initialize $\bm{\mathscr{q}}_{t_n} = \bm{q}\in \mathbb{R}$
\STATE Initialize $\iota_1 = 0$
\STATE Initialize $\iota_2 = 0$
\STATE Initialize $\iota_3 = 0$
\STATE Initialize drift $\partial U_{t_n}$ for $n\in\{0,1,\ldots,N\}$
\FOR{ $i$ in $n,\ldots, N-1$}
\STATE Sample Brownian noise: $\epsilon \sim \mathcal{N}(0,1)$
\STATE Compute the BEL weights:
 \IF{i == n}
\STATE $\iota_2 = \iota_2 +  \sqrt{\delta_t}\,\epsilon$
 \ELSE
\STATE $\iota_1 = \iota_1 +  \mu\,\delta_t \,(\partial^2 U_{t_n})(\bm{\mathscr{q}}_{t_n})$
\STATE $\iota_2 = \iota_2 + \sqrt{\delta_t}\,\epsilon\,e^{\,-\iota_1}$
 \ENDIF
\STATE $\iota_3 = \iota_3+ 
\delta_{t} \,\|(\partial_q U_{t_{i}} )(\bm{\mathscr{q}}_{t_{i}})\|^2\,\iota_2$ \\
\STATE Evolve one step of \eqref{od:FP}: $\bm{\mathscr{q}}_{t_{i+1}} = \bm{\mathscr{q}}_{t_{i}} - \mu \,(\partial_q U_{t_{i}})(\bm{\mathscr{q}}_{t_{i}})\,\delta_t + \sqrt{\delta_{t}\frac{2\,\mu}{\beta}}\, \epsilon$ 
\ENDFOR
\STATE Return $\partial V_{t_n}(\bm{q}) =\sqrt{\dfrac{\beta\,\mu}{2}}\,\biggl(\dfrac{\varphi(\bm{\mathscr{q}}_{\tf})}{\tf-t_n}\,\iota_2 +\dfrac{\beta\,\mu}{4} \,\iota_3\biggr) $
\end{algorithmic}
\end{algorithm}
 \begin{figure}
     \centering
    \includegraphics[width=\linewidth]{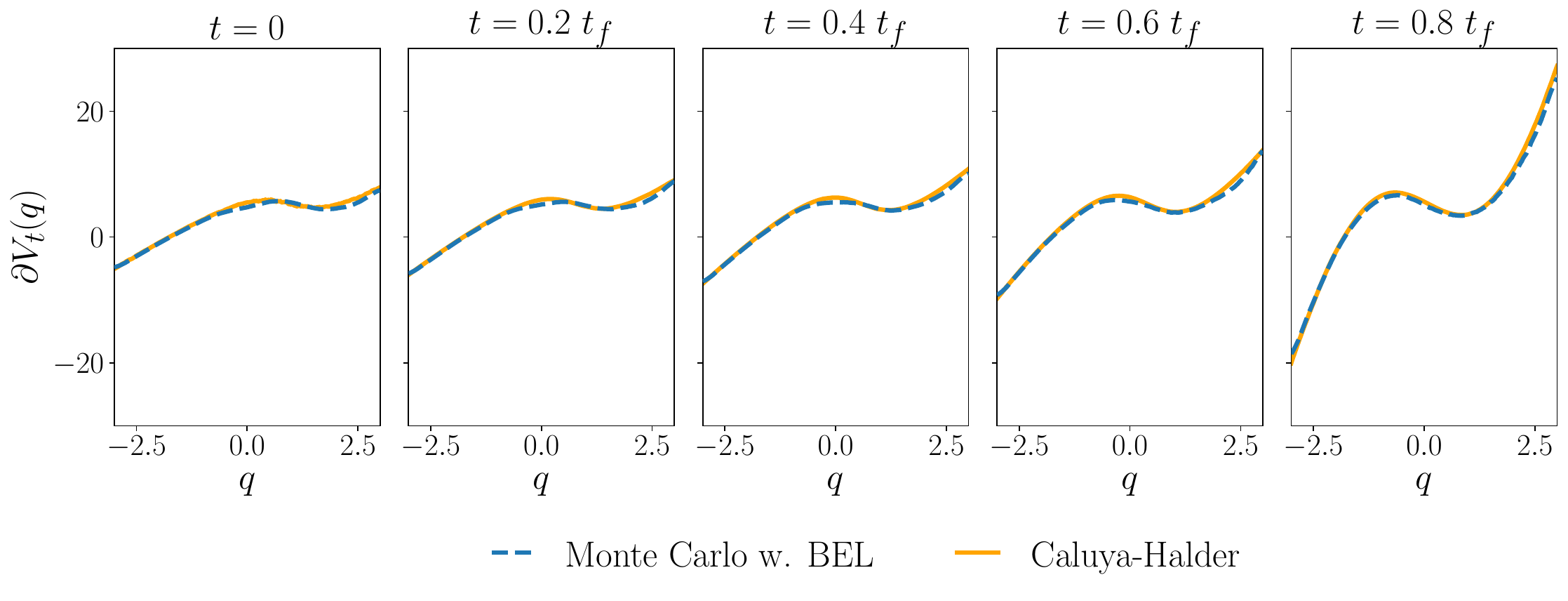}
     \caption{\label{fig:bel_dpe}Gradient of the value function, i.e. the gradient of the solution to the Hamilton-Jacobi-Bellman equation~\eqref{eq:hjb_coupled}, computed using the Bismut-Elworthy-Li formula (BEL) (dashed blue line) described in Algorithm \ref{alg_dvbel_over}. We sample  $10\ 000$ trajectories of the stochastic process \eqref{od:FP} from $500$ initial points in the interval $[-3,3]$, discretized by the Euler-Maruyama scheme with time step size $h=0.005$ and compute the BEL weights along the trajectories. The optimal control protocol $U$ and reference solution (orange) used is  computed by the iteration as in Fig.~\ref{fig:overdamped_fpk}. Numerical parameters and boundary conditions are as in Fig.~\ref{fig:overdamped_fpk}. We use $\ti = 0$, $\tf =0.2$ and $\mu=\beta=1$.}
 \end{figure}

\section{Application of Bismut-Elworthy-Li to degenerate diffusion}\label{sec:ud_bel}

For a degenerate diffusion we cannot directly apply (\ref{M:nd}) as it is because the expression is written in terms of the inverse of a degenerate matrix. Nevertheless the Bismut-Elworthy-Li formula continues to hold. To give an idea of how this comes about, we
consider the counterpart to (\ref{BEL:sde}), while referring to \cite{ZhaY2010} for the mathematically rigorous reader. Our starting point is
\begin{align}
\begin{split}
    	&	\mathrm{d}\bm{\mathscr{q}}_{t}= \bm{a}_{t}(\bm{\mathscr{x}}_{t})\,\mathrm{d}t
	\\
	&\mathrm{d}\bm{\mathscr{p}}_{t}=\bm{b}_{t}(\bm{\mathscr{x}}_{t})\,\mathrm{d}t+\mathsf{A}_{t}(\bm{\mathscr{x}}_{t})\,\mathrm{d}\bm{\mathscr{w}}_{t}
\end{split}
\label{dg:sde}
\end{align}
with 
\begin{align}
	\bm{\mathscr{x}}_{t}=\begin{bmatrix}
		\bm{\mathscr{q}}_{t}	\\  \bm{\mathscr{p}}_{t}
	\end{bmatrix}\colon \mathbb{R}_{+}\mapsto \mathbb{R}^{2\,d}
	\nonumber
\end{align}
and $\mathsf{A}$ a non-singular matrix field. The variational equation is
\begin{align}
	&	\mathrm{d}\bm{\mathscr{q}}_{t}^{\prime}= \bm{a}_{*t}(\bm{\mathscr{x}}_{t})\bm{\mathscr{x}}_{t}^{\prime}\,\mathrm{d}t
	\nonumber\\
	&\mathrm{d}\bm{\mathscr{p}}_{t}^{\prime}=\bm{b}_{*t}(\bm{\mathscr{x}}_{t})\bm{\mathscr{x}}_{t}^{\prime}\,\mathrm{d}t+\mathsf{A}_{*t}(\bm{\mathscr{x}}_{t})\bm{\mathscr{x}}_{t}^{\prime}\,\mathrm{d}\bm{\mathscr{w}}_{t}+\bm{\mathscr{h}}_{t}\,\mathrm{d}t
	\nonumber
\end{align}
where we suppose again
\begin{align}
	\bm{\mathscr{x}}_{s}^{\prime}=0
	\nonumber
\end{align}
We notice that we can always write the solution of the variational equation as
\begin{align}
\bm{\mathscr{x}}_{t}^{\prime}=\bm{\mathscr{x}}_{\star t,s}\bm{c}_{t}
	\nonumber
\end{align}
for some vector valued process $\bm{c}_{t}\colon \mathbb{R}\mapsto \mathbb{R}^{2\,d} $ such that
\begin{align}
	\bm{c}_{s}=0
	\nonumber
\end{align}
Upon differentiation, we readily verify that the self-consistency condition is
\begin{align}
	\bm{\mathscr{x}}_{\star \,t,s}\,\bm{\dot{c}}_{t}
	=\begin{bmatrix}
		0	\\  \bm{\mathscr{h}}_{t}
	\end{bmatrix}
	\nonumber
\end{align}
whose solution is
\begin{align}
	\bm{c}_{t}=\int_{s}^{t}\mathrm{d}u\,\bm{\mathscr{x}}_{\star u,s}^{-1}(\bm{x})\begin{bmatrix}
		0	\\  \bm{\mathscr{h}}_{u}
	\end{bmatrix}
	\label{dg:h}
\end{align}
We now avail us of the fact that the above relations hold for a sufficiently regular but otherwise arbitrary vector field $\bm{\mathscr{h}}_{u}$ and choose it such that
\begin{align}
	\bm{c}_{\tf}=\begin{bmatrix}
		0	\\ \bm{v}
	\end{bmatrix}
	\label{dg:bc}
\end{align}
Here $\bm{v}$ is a unit vector that specifies the direction of the gradient in Bismut-Elworthy-Li formula. Namely, given
\begin{align}
	V_{t}(\bm{x})=\operatorname{E}\left(\varphi(\bm{\mathscr{x}}_{\tf})\,\Big{|}\,\bm{\mathscr{x}}_{t}=\begin{bmatrix}
		\bm{q}	\\  \bm{p}
	\end{bmatrix}\right)
	\nonumber
\end{align}
the Bismut-Elworthy-Li formula \cite{ZhaY2010} continues to hold according to the chain of identities
\begin{align}
(\bm{v}\cdot
\bm{\partial}_{\bm{p}})V_{t}(\bm{x})
	& =\operatorname{E}\Bigl{\langle}\bm{\mathscr{x}}_{*\tf,t} \,\bm{c}_{\tf} \,, (\bm{\partial}\varphi)(\bm{\mathscr{x}}_{\tf})\Bigr{\rangle}
 =\operatorname{E}\Bigl{\langle} \bm{\mathscr{x}}_{\tf}^{\prime} \,, (\bm{\partial}\varphi)(\bm{\mathscr{x}}_{\tf})\Bigl{\rangle}
 \nonumber\\
& =\operatorname{E}\left(\varphi(\bm{\mathscr{x}}_{\tf})\int_{t}^{\tf}\left \langle\,\mathrm{d}\bm{\mathscr{w}}_{u}\,,\mathsf{A}_{u}^{-1}(\bm{\mathscr{x}}_{u})\bm{\mathscr{h}}_{u}\,\right\rangle\,\Big{|}\,\bm{\mathscr{x}}_{t}=\begin{bmatrix}
		\bm{q}	\\  \bm{p}
	\end{bmatrix}\right)
	\label{dg:BEL}
\end{align}
provided the conditions (\ref{dg:h}), (\ref{dg:bc}) are satisfied.

Similarly, we obtain a representation of the derivative with respect to the $\bm{q}$ variables by alternative choices of $\bm{\mathscr{h}}_{t}$
such that
\begin{align}
		\bm{c}_{\tf}=\begin{bmatrix}
		\bm{v} \\ 0
		\end{bmatrix}
	\nonumber
\end{align}

\subsection{A strategy for the explicit construction of a variational field enforcing the boundary conditions}
\label{sec:control}

Drawing from \cite{ZhaY2010}, we present a straightforward way to construct a variational field on the interval $[t,\tf]$ such that e.g. (\ref{dg:bc}) holds true.
Let
\begin{align}
	\mathds{H}_{u}:=\frac{\mathrm{d} \bm{\mathscr{x}}_{*\,u,t}}{\mathrm{d} u}\,\bm{\mathscr{x}}_{*\,u,t}^{-1}
	\nonumber
\end{align}
For clarity, we drop the subscripts $t,\tf$ in the following. However, there is still an implicit dependence on these parameters, with $u$ taking values in $[t,\tf]$. Consider the differential system
\begin{align}
	\begin{bmatrix}
		\bm{\dot{g}}_{u}	\\  \bm{\dot{f}}_{u}
	\end{bmatrix}
	=\mathds{H}_{u}
	\begin{bmatrix}
		\bm{g}_{u}	\\  \bm{\ell}_{u}
	\end{bmatrix}
	\label{dg:aux}
\end{align}
with $\bm{\ell}_{u}$ arbitrarily assigned (but sufficiently regular) and $\bm{g}_{u}$, and $\bm{f}_{u}$ determined by the identity. Then
\begin{align}
	\bm{\dot{\mathscr{x}}}_{u}^{\prime}+\begin{bmatrix}
		\bm{\dot{g}}_{u}	\\  \bm{\dot{f}}_{u}
	\end{bmatrix}
	=\mathds{H}_{u}\left(
	\bm{\mathscr{x}}_{u}^{\prime}+
	\begin{bmatrix}
		\bm{g}_{u}	\\  \bm{\ell}_{u}
	\end{bmatrix}
	\right)+\begin{bmatrix}
		0	\\  \bm{\mathscr{h}}_{u}
	\end{bmatrix}
	\nonumber
\end{align}
holds by construction for $u\in [t,\tf]$. Hence, if we require
\begin{align}
	\bm{\dot{\ell}}_{u}=\bm{\dot{f}}_{u}- \bm{\mathscr{h}}_{u}
	\nonumber
\end{align}
we see that
\begin{align}
	\bm{\mathscr{y}}_{u}=\bm{\mathscr{x}}_{u}^{\prime}+
	\begin{bmatrix}
		\bm{g}_{u}	\\  \bm{\ell}_{u}
	\end{bmatrix}
	\nonumber
\end{align}
satisfies
\begin{align}
&
\begin{split}
	&	\bm{\dot{\mathscr{y}}}_{u}=\mathds{H}_{u}\bm{\mathscr{y}}_{u}
	\\
	&\bm{\mathscr{y}}_{t}=\begin{bmatrix}
		\bm{g}_{t}	\\  \bm{\ell}_{t}
	\end{bmatrix}
	\end{split}
 && \forall\,u\in [t,\tf]
 \nonumber
\end{align}
 We solve (\ref{dg:aux}) with the {\textquotedblleft}initial condition{\textquotedblright}
\begin{align}
\bm{g}_{t}=0
	\nonumber
\end{align}
As a consequence, we arrive at
\begin{align}
	\bm{\mathscr{x}}_{u}^{\prime}=\bm{\mathscr{x}}_{*u,t}(\bm{x})\begin{bmatrix}
		0	\\  \bm{\ell}_{t}
	\end{bmatrix}
	-\begin{bmatrix}
		\bm{g}_{u}	\\  \bm{\ell}_{u}
	\end{bmatrix}
	\nonumber
\end{align}
The identity we just obtained shows that in order to get a representation of the gradient according to (\ref{dg:BEL}), we must restrict the choice of vector fields $\bm{\ell}_{u}$ to those satisfying the boundary conditions
\begin{align}
&   \bm{\ell}_{t}=\bm{v} \qquad \& \qquad \bm{\ell}_{\tf}=0
\end{align}
so that at time $u=\tf$ 
\begin{align}
	\bm{\mathscr{x}}_{\tf}^{\prime}=\bm{\mathscr{x}}_{*\tf,t}(\bm{x})\begin{bmatrix}
		0	\\  \bm{v}
	\end{bmatrix}
	\label{dg:wanted}
\end{align}
holds true. Once all the above conditions are satisfied, we can determine the right hand side of (\ref{dg:BEL}) from
\begin{align}
&    \bm{\mathscr{h}}_{u}=\bm{\dot{f}}_{u}-\bm{\dot{\ell}}_{u} && \forall u \,\in \,[t,\tf]
    \label{dg:explicit}
\end{align}

\subsubsection{A case of particular interest}

There are an infinite number of ways to choose $\bm{\ell}_{u}$ such that the condition (\ref{dg:wanted}) holds true.
We detail here a choice of particular interest for physics. Let us consider the generator of the linearized dynamics around a path solution of (\ref{ud:sde})
\begin{align}
	\mathds{H}_{u}(\bm{\mathscr{x}}_{u})=\begin{bmatrix}
		0 & \frac{\operatorname{1}}{m}	\\[0.3cm]  -(\bm{\partial}\otimes\bm{\partial} U_{u})(\bm{\mathscr{q}}_{u}) & -\frac{\operatorname{1}}{\tau}
	\end{bmatrix}
	\nonumber
 \end{align}
In such a case,  the instantiation of (\ref{dg:aux}) is the differential system
\begin{equation}
	\begin{split}
	&	\bm{\dot{g}}_{u}=\frac{1}{m}\bm{\ell}_{u}
	\\
	&    \bm{\dot{f}}_{u}=-(\bm{\partial}\otimes\bm{\partial} U_{s})(\bm{\mathscr{q}}_{u})\,\bm{g}_{u}-\frac{1}{\tau} \bm{\ell}_{u}
	\\
	& \bm{\dot{f}}_{u}-\bm{\mathscr{h}}_{u}=\bm{\dot{\ell}}_{u}
\end{split}
\nonumber
\end{equation}
We assign
\begin{align}
	\bm{\ell}_{u}=\bm{v}\,\frac{\tf -u}{\tf-t}+\bm{v}_{\mathfrak{1}}\,\frac{(\tf -u)(u-t)}{2\,(\tf-t)^{2}}
	\nonumber
\end{align}
and obtain
\begin{align}
	\bm{g}_{u}=\bm{v}\int_{t}^{u}\mathrm{d}s\,\dfrac{\tf -s}{m\,(\tf-t)}+\bm{v}_{\mathfrak{1}}\int_{t}^{u}\mathrm{d}s\,\dfrac{(\tf -s)(s-t)}{2\,m\,(\tf-t)^{2}}
	\nonumber
\end{align}
We fix $\bm{v}_{\mathfrak{1}} $ by requiring
\begin{align}
	\bm{g}_{\tf}=0\Longrightarrow \bm{v}_{\mathfrak{1}}=-6\,\bm{v}
	\nonumber
\end{align}
We thus obtain
\begin{align}
	&\bm{\ell}_{u}=\bm{v}\,\dfrac{(\tf-u)(\tf+2\,t-3\,u)}{(\tf-t)^{2}}
	\\	
	&	\bm{g}_{u}=\bm{v}\,\dfrac{(\tf-u)^{2}(u-t)}{m\,(\tf-t)^{2}}
 \nonumber
\end{align}
and therefore
\begin{align}
	\bm{\dot{f}}_{u}=-\frac{1}{\tau}\bm{v}\,\dfrac{(\tf-u)(\tf+2\,t-3u)}{(\tf-t)^{2}}- (\bm{\partial}\otimes\bm{\partial} U_{u})(\bm{\mathscr{q}}_{u})\,\bm{v}\,\dfrac{(\tf-u)^{2}(u-t)}{m\,(\tf-t)^{2}}
	\nonumber
\end{align}

\subsection{Analytical example}\label{sec:an_hjb_ud}

We return to the elementary case (\ref{exa:sde}), (\ref{exa:pde}) but set $\eta=0$. The momentum gradient (\ref{exa:BEL}) does not depend on $\eta$, yet the application of the Bismut-Elworthy-Li formula requires the inverse of the volatility which in turn appears to depend on $\eta$. 
As in the present example, the potential is identically vanishing 
\begin{align}
U=0
\end{align}
we arrive at
\begin{align}
	&	\partial_{p}V_{t}(p,q)=
	\nonumber\\
	&
\operatorname{E}\left(
	\biggl(\mathscr{p}_{t}e^{-\frac{\tf-t}{\tau}}+\sqrt{\frac{2\,m}{\beta\,\tau}}\int_{t}^{\tf}\mathrm{d}w_{s}^{(\mathfrak{2})}\,e^{-\frac{\tf-s}{\tau}}\biggr)
	\int_{t}^{\tf}
	\sqrt{\frac{\beta\,\tau}{2\,m}}\mathrm{d}w_{u}^{(\mathfrak{2})} \left(-\frac{\mathscr{\ell}_{u}}{\tau}-\dot{\mathscr{\ell}}_{u}
	\right)
	\,\bigg{|}\,\mathscr{q}_{t}=q,\mathscr{p}_{t}=p\right)
 \nonumber
\end{align}
Using again the properties of stochastic integrals, the expectation value reduces to
\begin{align}
\partial_{p}V_{t}(p,q)=-\int_{t}^{\tf}\mathrm{d}s\,\frac{\mathrm{d}}{\mathrm{d} s}\left(e^{-\frac{\tf-s}{\tau}}\mathscr{\ell}_{s}\right)
	=-\mathscr{\ell}_{\tf}+e^{-\frac{\tf-t}{\tau}}\mathscr{\ell}_{t}
\nonumber
\end{align}
whence we recover the correct expression of the gradient once we recall the boundary conditions imposed on the function $\mathscr{\ell}$ in $[t,\tf]$.

This example also indicates that the existence of the Bismut-Elworthy-Li formula for Langevin-Kramers equations of the form (\ref{dg:sde}) can be recovered from the limit $\eta$ tending to zero of a non-degenerate model owing to the vanishing of products of It\^o stochastic integrals with respect to independent Wiener processes. 

\subsection{Numerical Example}\label{sec:num_hjb_ud}

We demonstrate here a numerical example of using the Bismut-Elworthy-Li formula to find the gradient of a value function satisfying the  Hamilton-Jacobi-Bellman equation \eqref{num_ud_hjb}. We look at the case where the initial and final conditions assigned on the density are Gaussian distributions. In the case of Gaussian boundary conditions, we can determine the value function and optimal control protocol in the underdamped dynamics as the numerical solution of a system of differential equations; see Section IV in \cite{SaBaMG2024}. The value function $V_t$ is quadratic in the momentum and position variables, and in the two-dimensional phase space case reads
\begin{align*}
    V_t(p,q) = v^{(0)}_t + v_t^{(p)} p + v_t^{(q)} q + \dfrac{1}{2}\biggl(v_t^{(p,p)}p^2+ 2 v_t^{(p,q)}pq+v_t^{(q,q) }q^2\biggr)
\end{align*}
for time-dependent coefficients $v^{(0)}_t,\,v_t^{(p)} , \,v_t^{(q)} ,\,v_t^{(p,p)}, \,v_t^{(p,q)}$ and $v_t^{(q,q)}$ found as in Section IV of \cite{SaBaMG2024}. The solution of \eqref{BEL:dp} can be found as
\begin{align*}
    V_t(\bm{x}) = \operatorname{E}\Bigl(\varphi(\bm{\mathscr{x}}_{\tf}) + \dfrac{\beta\,\tau}{4\,m} \int_t^{\tf} \dd s \, \|\partial_q U_s(\bm{\mathscr{q}}_s)\|^2\,\big{|}\,\bm{\mathscr{x}}_t = \bm{x}\Bigr)
\end{align*}
Applying the Bismut-Elworthy-Li formula with $\bm{\mathscr{h}}_{t}$ gives the following expression for the gradient of the value function with respect to momentum
\begin{align*}
    \partial_p V_t(\bm{x}) &= \sqrt{\dfrac{\beta\,\tau}{2\,m}}\,\operatorname{E}\biggl(\varphi(\bm{\mathscr{x}}_{\tf})\int_t^{\tf} \langle\dd \bm{\mathscr{w}}_s \,,\, \bm{\mathscr{h}}_{s}\rangle  \\
    &\qquad+ \dfrac{\beta\,\tau}{4\,m} \int_t^{\tf} \dd s\,\|\partial U_s(\bm{\mathscr{q}}_s)\|^2\int_t^s\langle\dd \bm{\mathscr{w}}_u \,,\, \bm{\mathscr{h}}_{u} \rangle \,\Big{|}\, \bm{\mathscr{x}}_t = \bm{x}\biggr) 
\end{align*}
where
\begin{align}\label{num_UD_BEL:hu}
\bm{\mathscr{h}}_{u} &=-\bm{\dot{\mathscr{\ell}}}_{u}-\frac{\bm{\mathscr{\ell_{u}}}}{\tau}-(\bm{\partial}\otimes \bm{\partial} U_{u})(\bm{\mathscr{q}}_{u} )\,\bm{\mathscr{g}}_{u}
    \nonumber\\
& =-\bm{v} \,\dfrac{1}{\tau (\tf-t)^2}\,\biggl((\tf-u)(\tf+2t-3u)-\tau\,(4\tf-6u+2t) \\
&\qquad\qquad +  \frac{\tau}{m}\,(\bm{\partial}\otimes \bm{\partial} U_{u})(\bm{\mathscr{q}}_{u})(\tf-u)^{2}(u-t)\biggr)
\nonumber
\end{align}
The computation is summarized in Algorithm~\ref{alg_dvbel_under}.
\begin{figure}
    \centering
    \includegraphics[width=0.5\linewidth]{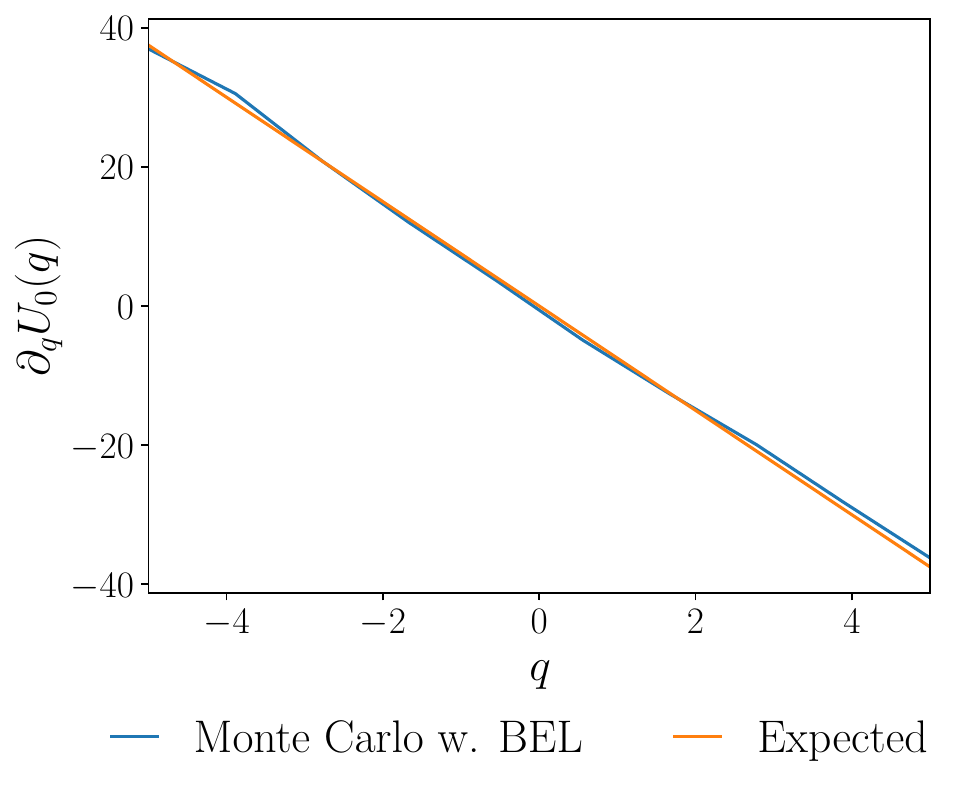}
    \caption{The gradient of the optimal control potential minimizing the Kullback-Leibler divergence~\eqref{eq:kl_divergence} in the underdamped dynamics. We compute the stationarity condition~\eqref{num_ud_stationarity_condition} at $t=0$ using the gradient of the solution of the Hamilton-Jacobi-Bellman equation~\eqref{BEL:dp} using the Bismut-Elworthy-Li formula (Monte Carlo w. BEL) (blue line) in Algorithm~\ref{alg_dvbel_under}. The optimal control protocol $U$ and terminal condition $\varphi$ of~\eqref{BEL:dp} are found using numerical integration of the system of equations described in Section IV of~\cite{SaBaMG2024}, using a fourth order co-location method from the DifferentialEquations.jl library~\cite{differentialequationsjl}. We use Gaussian boundary conditions: the initial and final position and momentum means are set as zero; the initial and final cross-correlation is zero; the initial variances are set to $1$; the final position variance is $1.7$, and the final momentum variance is $1$. We sample $10\ 000$ independent trajectories of the stochastic process~\eqref{ud:sde} started from $500$ sample points in the interval $[-5,5]\times[-5,5]$ using an Euler-Maruyama discretization with time-step $h=0.01$. We use $\ti = 0$, $\tf =1$ and $\beta=\tau = m = 1$.}
    \label{fig:bel_underdamped}
\end{figure}

\begin{algorithm}
\caption{BEL for degenerate diffusion}
\label{alg_dvbel_under}
\begin{algorithmic} 
\STATE Initialize $\bm{\mathscr{q}}_{t_n} = q\in \mathbb{R}$
\STATE Initialize $\bm{\mathscr{p}}_{t_n} = p\in \mathbb{R}$\STATE Initialize $\iota_1 = 0$
\STATE Initialize drift function $\partial_q U_{t}$
\STATE Initialize $h(u,t,T)$ as in \eqref{num_UD_BEL:hu}
\FOR{ $i$ in $n,\ldots, N-1$}
\STATE Sample Brownian noise: $\epsilon_{i} \sim \mathcal{N}(0,\delta_t)$
\IF{$i>n$}
\STATE Add to running cost: $\iota_1 = \iota_1 + \delta_t\|\partial_q U_{t_{i}} (\bm{\mathscr{q}}_{t_{i}})\|^2\Bigl(\sum_{j=n}^{i} \epsilon_j\,h(t_j,t_n,t_i) \Bigr)$
\ENDIF
\STATE Evolve one step of \eqref{ud:sde}: $\begin{cases}
    
\bm{\mathscr{q}}_{t_{i+1}} = \bm{\mathscr{q}}_{t_{i}} + \frac{\bm{\mathscr{p}}_{t_{i}}}{m}\,\delta_t \\[5pt]
\bm{\mathscr{p}}_{t_{i+1}} = \bm{\mathscr{p}}_{t_{i}} - \Bigl(\frac{\bm{\mathscr{p}}_{t_{i}}}{\tau} + \partial_q U_{t_{i}}(\bm{\mathscr{q}}_{t_{i}})\Bigr)\,\delta_t + \sqrt{\frac{2\,m}{\tau\beta}}\, \epsilon_i\end{cases}$ 
\ENDFOR
\STATE Return $\partial_p V_{t_n}(\bm{\mathscr{x}}) =\sqrt{\dfrac{\beta}{2}}\biggl(\varphi(\bm{\mathscr{x}}_{\tf})\Bigl(\sum_{j=n}^{N-1} \epsilon_{j} \, h(t_j,t_n,\tf)\Bigr) +\dfrac{\beta}{4} \,\iota_1 \biggr)$
\end{algorithmic}
\end{algorithm}

\section{Application to Machine Learning}
In this section, we return to the overdamped dynamics and demonstrate an application of numerical methods we discuss above. We present a prototype example for the optimal control problem in the overdamped dynamics of minimizing the Kullback Leibler divergence \eqref{eq:kl_divergence}. Inspired by the seminal works \cite{EHaJe2017,HaJeE2018}, we model the optimal control protocol by a neural network, and use gradient descent to iteratively update it based on the stationarity condition \eqref{stationarity_condition}. 

As before, we formulate the problem in terms of a Bismut-Pontryagin cost functional. Additionally, we enforce the assigned boundary conditions (initial and final conditions on the density of the form \eqref{eq:boundary_conditions}) through a Lagrangian multiplier $\lambda$. This gives 
\begin{align}\label{eq:bp_cost}
 	\mathcal{A}[\mathtt{p},U,V]&=\int_{\mathbb{R}^{d}}\mathrm{d}^{d}\bm{q}\left(- V_{\tf}(\bm{q})\mathtt{p}_{\tf}(\bm{q})+\lambda(\bm{q})\,\big{(}\mathtt{p}_{\tf}(\bm{q})-P_f(\bm{q})\big{)}
 	+V_{\ti}(\bm{q})P_{\iota}(\bm{q})\right)
 	\nonumber\\
 	&	+
 	\int_{\ti}^{\tf}\mathrm{d}t\,\int_{\mathbb{R}^{d}}\mathrm{d}^{d}\bm{q}\,\mathtt{p}_{t}(\bm{q})\biggl(\frac{\beta\mu}{4}\|(\bm{\partial}U_{t})(\bm{q})\|^{2}+\Bigl(\partial_{t}-\mu\langle(\bm{\partial}_q U)(\bm{q})\,,\,\bm{\partial}_q\rangle + \dfrac{\mu}{\beta}\bm{\partial}_q^2\Bigr)V_{t}(\bm{q})\biggr)
 	\nonumber
 \end{align}
Taking stationary variation with respect to the density $\mathrm{p}$, control protocol $U$ and value function $V$ yields the coupled partial differential equations \eqref{eq:fpk_coupled} and \eqref{eq:hjb_coupled}, and the stationarity condition \eqref{stationarity_condition}. We identify 
\begin{equation}
    V_{\tf}(\bm{q}) = \lambda(\bm{q})
\end{equation}
and the following update rule
       \begin{equation}
\label{num_value_update_rule}\lambda^{(\text{new})}(\bm{q}) = \lambda^{(\text{old})}(\bm{q}) + \gamma_1 \,\biggl(\log{\frac{\mathrm{p}_{\tf}(\bm{q})}{P_f(\bm{q})}}\biggr)
\end{equation}
chosen in this way to preserve the integrability conditions of the value function. 
The stationarity condition gives an update rule for the drift of the control protocol
\begin{equation}
\label{num_drift_update_rule}(\bm{\partial} U_t)^{(\text{new})}(\bm{q}) = (\bm{\partial}U_t)^{(\text{old})}(\bm{q}) - \gamma_2 \, \biggl((\bm{\partial}U_t)^{(\text{old})}(\bm{q})  - \frac{\beta}{2}(\bm{\partial} V_t)(\bm{q})\biggr)
\end{equation}
The parameters $\gamma_1,\gamma_2>0$ control the step sizes of the gradient descent, known as a learning rate. The update for the Lagrange multiplier is a gradient ascent rather than descent \cite{constrained_diff_opt}. 

The right hand sides of both \eqref{num_value_update_rule} and \eqref{num_drift_update_rule} can be computed using Monte Carlo integration techniques discussed in this note. With appropriate parametrization of the gradient of the control protocol and the Lagrangian multiplier $\lambda$, the method could therefore scale to high dimensions. In this prototype example, we use a polynomial regression for fitting $\lambda$ and a neural network for the gradient of the control protocol. The polynomial regression could be replaced with any suitable parametrization, in particular, with a second neural network. 

The gradient of the optimal control protocol $\bm{\partial} U_t$ is modelled by a neural network, denoted by $\mathbb{U}_t$. We use a feed-forward neural network: connected layers, representing affine transformations with non-linear functions (known as activation functions) between them. The neural network has a set of parameters (weights and biases) associated with the layers, which we denote by $\Theta$. The network takes the time $t$ and space coordinates $\bm{q}$ as input. Using a neural network allows for evaluating the optimal control protocol on new space coordinates without using interpolation, meaning that it can easily be used as the drift in the computation of the density and value functions using Algorithms~\ref{alg:fpkgirsanov} and~\ref{alg_dvbel_over}. 

The training process can be summarized as follows. Firstly, the functions $\lambda$ and $\mathbb{U}$ with a set of parameters $\Theta^{(0)}$ are initialized. We use these to find the final density $\mathrm{p}_{\tf}$ under this drift with Algorithm~\ref{alg:fpkgirsanov}. The Lagrange multiplier is updated using \eqref{num_value_update_rule}. The new $\lambda$ is used as the terminal condition of the value function. We then use Algorithm~\ref{alg_dvbel_over} to compute the gradient of the value function, $\partial_q V_t$, using the current drift and terminal condition. The neural network parameters $\Theta$ are updated so that the new drift satisfies \eqref{num_drift_update_rule}. Under the updated drift, the final density is recomputed and the process is repeated until convergence. The whole iteration is summarized in Algorithm~\ref{alg2:all}.

\begin{figure}
     \centering
    \includegraphics[width=\linewidth]{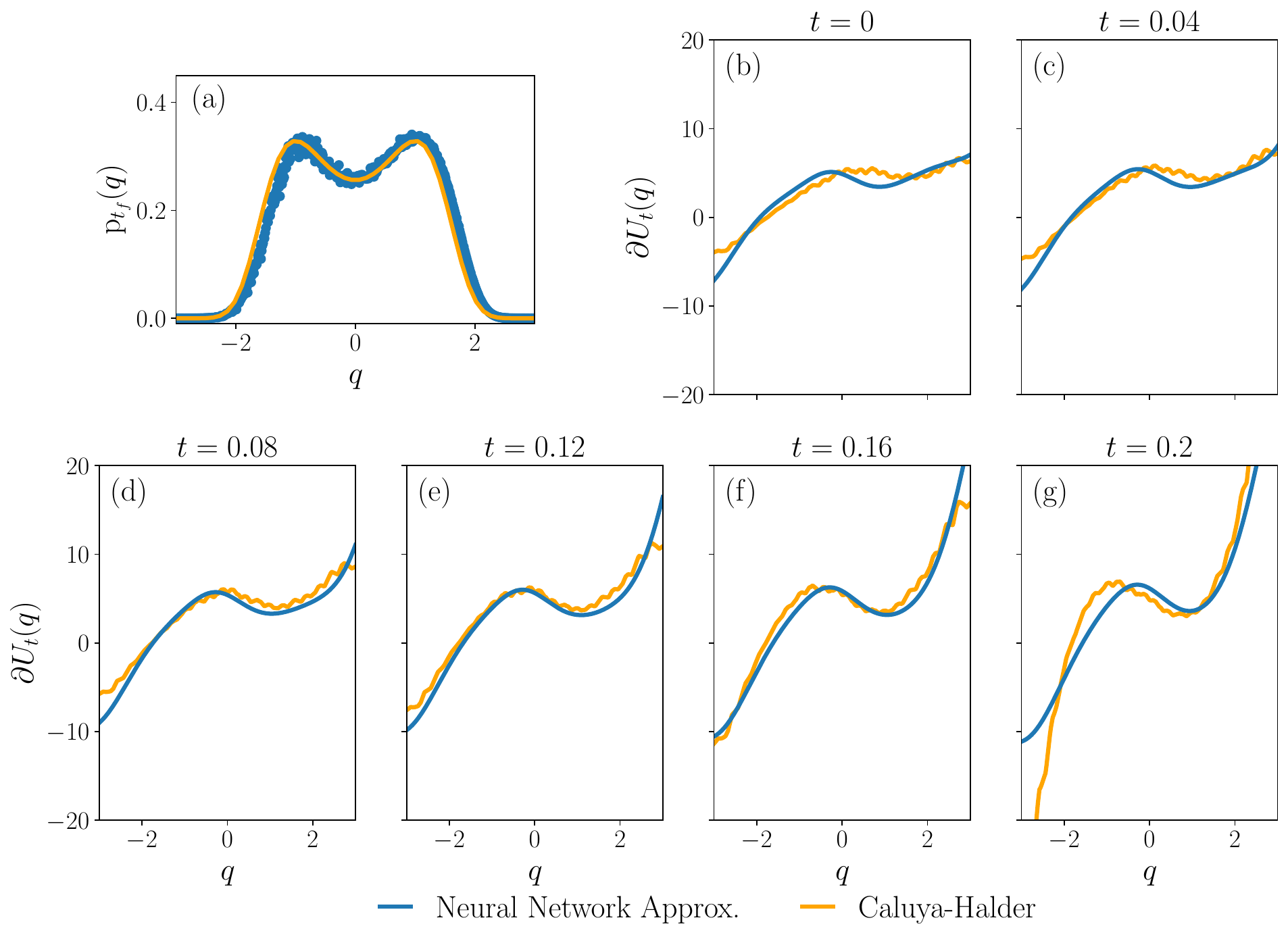}     \caption{\label{fig:nnresults} Solution of the optimal control problem minimizing the Kullback-Leibler divergence from a free diffusion in a fixed time interval in the overdamped dynamics. The gradient of the control protocol is parameterized by a neural network and trained using the iteration described in Algorithm~\ref{alg2:all}. Panel~(a) shows the final boundary condition obtained by integrating the Fokker-Planck equation \eqref{eq:fpk_coupled} using the trained neural network as the drift in Algorithm~\ref{alg:fpkgirsanov} (blue) against the assigned final boundary condition (orange). Panels~(b)-(g) show the output of the neural network (blue) after training to estimate the gradient of the optimal control protocol against a reference solution~\cite{CaHa2022} (orange). \newline
    We use the assigned boundary conditions as in Fig.~\ref{fig:overdamped_fpk}, with $\beta=\mu=1$, $\ti=0$ and $\tf=0.2$. The gradient optimal control protocol is parametrized by a fully connected feed-forward neural network with one input layer of four neurons, one hidden layer of ten neurons and one output layer. The swish ($x \longmapsto x\sigma(x)$) activation function is used between the input and hidden, and hidden and output layers. Weights and biases are initialized using Glorot normal initialization, and Glorot uniform initialization for the output layer. The Lagrange multiplier $\lambda$ function is approximated by fitting a polynomial of degree $6$ and initialized with all coefficients set to $0$. At each iteration, $512$ points are sampled uniformly from the interval $[-3,3]$. The gradient of the value function is computed using Alg.~\ref{alg_dvbel_over} with $\mathbb{U}$ as the drift. Each computation uses $10$ independent simulated trajectories of the associated SDE using an Euler-Maruyama discretization and time-step $0.005$. The final density is computed using Alg.~\ref{alg:fpkgirsanov} with $\mathbb{U}$ as the drift. Each computation uses $100$ independent Monte Carlo trajectories from each sample point using an Euler-Maruyama discretization and time-step $0.005$. \newline
    The neural network $\mathbb{U}$ is trained in four phases as follows. The first phase is $20$ full iterations of Alg.~\ref{alg2:all} with $100$ updates to the parameters $\Theta$ per iteration using stochastic gradient descent with learning rate $\gamma_2 = 10^{-3}$. At each iteration, the Lagrange multiplier $\lambda$ is recomputed using \eqref{num_value_update_rule} with $\gamma_1= 0.1$. In the second phase, we make $20$ full iterations with $100$ updates to $\Theta$ according to \eqref{num_drift_update_rule} using stochastic gradient descent and learning rate $\gamma_2 = 10^{-4}$ per iteration. The Lagrange multiplier is recomputed once at each iteration using $\gamma_1 = 10^{-2}$. In the third phase, we make $20$ full iterations with $400$ updates to $\Theta$ using stochastic gradient descent with learning rate $\gamma_2 = 10^{-5}$ per iteration. In the fourth phase, we make $20$ full iterations with $400$ update steps to $\Theta$ per iteration using the ADAM \cite{kingma2014} optimizer and $\gamma_2= 10^{-4}$. The code is written in the Julia programming language, using especially the Flux.jl~\cite{fluxjl2,Fluxjl} and Polynomials.jl libraries.}
 \end{figure}

The results of Algorithm~\ref{alg2:all} are illustrated in Fig.~\ref{fig:nnresults}. We show the final density in panel (a). Panels (b)-(g) show the approximation of the gradient of the control protocol by the trained neural network.

\begin{algorithm}
\caption{Learning an Optimal Control Protocol by Gradient Descent}
\label{alg2:all}
\begin{algorithmic} 
\STATE Initialize a neural network $\mathbb{U}$ with parameters $\Theta^{(0)}$
\STATE Initialize $\lambda^{(0)}$ as a polynomial with coefficients zero
\STATE Initialize $\gamma_1,\gamma_2$ learning rates
\WHILE{$\ell \leq L$ max iters}
\STATE Initialize a batch $\bm{q} = \{q_k\}_k$ of $K$ points
\STATE Compute $\mathrm{p}^{(\ell)}_{\tf}(\bm{q})$ using Algorithm~\ref{alg:fpkgirsanov} with $\mathbb{U}$ as drift
\STATE Update $\lambda^{(\ell +1)}(\bm{q}) = \lambda^{(\ell)}(\bm{q}) + \gamma_1\biggl(\log\frac{\mathrm{p}^{(\ell)}_{\tf}(\bm{q})}{P_f(\bm{q})}\biggr)$ 
\STATE Set $V^{(\ell+1)}_{\tf}(\bm{q}) = \lambda^{(\ell+1)}(\bm{q})$
\FOR{$t_n$, $n=0,\ldots,N$}
\STATE Compute $(\partial V^{(\ell+1)}_{t_n})(\bm{q})$ using Algorithm~\ref{alg_dvbel_over} with $\mathbb{U}$ as drift
\STATE Update $\Theta^{(\ell+1)}$ such that $\biggl{\|}\mathbb{U}(\bm{q},t_n;\,\Theta^{(\ell+1)}) - \dfrac{\beta}{2}\,(\partial V^{(\ell+1)}_{t_n})(\bm{q})\biggr{\|}^2$ is minimized
\ENDFOR
\ENDWHILE
\RETURN  Approximation of the gradient of the optimal control protocol $\mathbb{U}$
\end{algorithmic}
\end{algorithm}

\section{Conclusions}

In this note, we discuss two integration methods for partial differential equations which frequently appear in optimal control problems. We show how we can use the Girsanov theorem such that a Fokker-Planck equation driven by a mechanical potential can be integrated by taking a numerical expectation of Monte Carlo trajectories of an auxiliary stochastic process. This method can be applied when the auxiliary stochastic process is non-degenerate or degenerate. Secondly, we use the Bismut-Elworthy-Li formula to find expressions for the gradient of the value function satisfying a Hamilton-Jacobi-Bellman equation. We show this for both a non-degenerate and degenerate diffusion.

The discussed numerical methods are supported by computational examples. We examine the dynamic Schr\"odinger bridge problem, or the minimization of the Kullback-Leibler divergence from a free diffusion while satisfying boundary conditions on the density at the initial and final time. For the overdamped dynamics, our integration shows good agreement with the iterative approach of Caluya and Halder \cite{CaHa2022} in Figs.~\ref{fig:overdamped_fpk} and \ref{fig:bel_dpe}. In the underdamped case, we integrate the associated Fokker-Planck equation to support the consistency of the multiscale perturbative approach used in  \cite{SaBaMG2024}. In particular, we compute an estimate of the evolution of the joint density function of the system state for this problem in Fig.~\ref{fig:fpk_underdamped}. We also verify the stationarity condition using the Bismut-Elworthy-Li for a degenerate diffusion in Fig.~\ref{fig:bel_underdamped}. Finally, we demonstrate an application of both integrations in a simple machine learning model in Fig.~\ref{fig:nnresults}. 

The optimal control problem discussed here has many applications. One possibility is application in machine learning, for instance in the development of diffusion models for image generation \cite{DeBoThHeAr2021}. Here, we find an optimal steering protocol between a noise distribution (e.g. a Gaussian) and a target (e.g. an image) by minimizing the Kullback-Leibler divergence. Optimal control problems in the underdamped dynamics are particularly interesting. Underdamped dynamics take into account random thermal fluctuations, noise and the effects of inertia, hence they are well suited to model non-equilibrium transitions at nanoscale. Models of certain biological systems require considering complex dynamics, for example, because of random external noise from the environment \cite{MaReOp2020}. Such models then result in non-linear partial differential equations, making them difficult to integrate. While the implementation of machine learning to solve an optimal control problem we use here is a prototype, it may be possible to extend it to a more general setting. Specifically, we have in mind transitions obeying underdamped dynamics and occurring at minimum entropy production such as those considered in \cite{SaBaMG2024_2}.

\paragraph{Acknowledgements} JS was supported by a University of Helsinki funded doctoral researcher position, Doctoral Programme in Mathematics and Statistics. The authors gratefully acknowledge Paolo Andrea Erdman for many useful discussions on topics related to the present work.

\paragraph{Data Availability}The codes used to generate numerical examples shown in this article are available in Github at \url{https://github.com/julia-sand/kldivergence}.

\printbibliography

\end{document}